\documentclass[12pt,thmsb]{article}
\usepackage{amsfonts}

\usepackage{graphicx}
\usepackage{sw20lart}

\input{tcilatex}

\begin{document}

\title{}

\begin{center}
{\LARGE A Model of Blood Flow in a Circulation Network}

\bigskip

\textbf{Weihua Ruan}$^{\dagger }$\textbf{, M.E. Clark}$^{\ddagger }$\textbf{%
, Meide Zhao}$^{\ddagger }$\textbf{\ and Anthony Curcio}$^{\ddagger }$

\medskip

$^{\dagger }$\textbf{Department of Mathematics, Computer Science and
Statistics,}

\textbf{Purdue University Calumet}

\smallskip

\textbf{and}

$^{\ddagger }$\textbf{VasSol, Inc.}
\end{center}

\smallskip \medskip

\paragraph{Abstract.}

We study a mathematical model of a blood circulation network which is a
generalization of the coronary model proposed by Smith, Pullan and Hunter.
We prove the existence and uniqueness of the solution to the
initial-boundary value problem and discuss the continuity of dependence of
the solution and its derivatives on initial, boundary and forcing functions
and their derivatives.

\smallskip

\section{\protect\smallskip Introduction\label{Int}}

In a recent paper \cite{SPH02}, Smith, Pullan and Hunter propose a
mathematical model of blood circulation in the coronary network,
and conduct a numerical analysis. In their model, major vessels
(with cross-sectional areas larger than a certain value) are
treated as a connected one-dimensional network, and small vessels,
such as arterioles, capillaries and venules, are treated as lumped
elements which are connected to the network of vessels. The flow
on vessels are assumed to be incompressible, homogeneous,
Newtonian, and has a small Reynolds number. Thus, the mass balance
equation and Navier-Stokes equation can be written to describe the
pressure and the flow rate on vessels. Equations on lumped
elements are written in analogy with the current and voltage in an
electric circuit. The result is an initial-boundary value problem
of a system of hyperbolic type partial differential equations
coupled at junctions of the network. Although the result of the
numerical analysis conducted in \cite{SPH02} matches closely with
measured data, the well-posedness problem of the system of partial
differential equations, that is, the existence, uniqueness and the
continuous dependence on initial and boundary data of the
solution, has not been established before. The main objective of
this paper is to establish the well-posedness. We prove that the
system is well-posed under certain natural conditions. This work
is an extension of our earlier work \cite{RCZC} on a model of
blood circulation in the brain. The main differences between the
two models are that the network configuration in \cite{RCZC} is
more complicated owing to the presence of Willis loops, but the
coupling junction conditions in the model of \cite{SPH02} are more
complicated due to the different formulation and the inclusion of
the capillaries and veinal system. We combine both features in a
more general system with the hope that our result will be useful
in the modelling of circulation systems of higher complexity,
including the whole body circulation system.

Before stating our system, let us briefly describe the model in \cite{SPH02}%
. Let $P_{i}$ and $R_{i}$ represent the pressure and radius on the $i$-th
vessel, respectively, and let $V_{i}$ be the cross-sectional average of the
axial component $v_{i,x}$ of the velocity on the $i$-th vessel. Assuming
that the radial component $v_{i,r}$ of the velocity is small compared to the
axial component $v_{i,x}$ of the velocity, one can write equations of mass
balance
\begin{equation}
\frac{\partial R_{i}}{\partial t}+V_{i}\frac{\partial R_{i}}{\partial x}+%
\frac{R_{i}}{2}\frac{\partial V_{i}}{\partial x}=0,  \label{massbal}
\end{equation}
and momentum balance
\begin{equation}
\frac{\partial V_{i}}{\partial t}+2\left( 1-\alpha _{i}\right) \frac{V_{i}}{%
R_{i}}\frac{\partial R_{i}}{\partial t}+\alpha _{i}V_{i}\frac{\partial V_{i}%
}{\partial x}+\frac{1}{\rho }\frac{\partial P_{i}}{\partial x}=\frac{2\nu }{%
R_{i}}\left[ \frac{\partial v_{i,x}}{\partial r}\right] _{r=R_{i}}.
\label{mombal}
\end{equation}
Here $\nu $ is the viscosity constant and
\[
\alpha _{i}=\frac{1}{R_{i}^{2}V_{i}^{2}}\int_{0}^{R_{i}}2rv_{i,x}^{2}dr
\]
is the energy quantity. Taking into consideration of no-slip boundary
condition ($v_{i,x}=0$ if $r=R_{i}$), the viscous axisymmetry ($\partial
v_{i,x}/\partial r=0$ if $r=0$), and the fact that $V_{i}$ is the
cross-sectional average of $v_{i,x}$, Smith, Pullan and Hunter propose the
velocity profile
\[
v_{i,x}\left( r,x\right) =\frac{\gamma _{i}+2}{\gamma _{i}}V_{i}\left(
x\right) \left[ 1-\left( \frac{r}{R_{i}}\right) ^{\gamma _{i}}\right]
\]
where $\gamma _{i}$ is a positive number. Using this profile and
the mass balance condition (\ref{massbal}), Eq. (\ref{mombal})
becomes
\begin{equation}
\frac{\partial V_{i}}{\partial t}+\left( 2\alpha _{i}-1\right) V_{i}\frac{%
\partial V_{i}}{\partial x}+2\left( \alpha _{i}-1\right) \frac{V_{i}^{2}}{%
R_{i}}\frac{\partial R_{i}}{\partial x}+\frac{1}{\rho }\frac{\partial P_{i}}{%
\partial x}=-\frac{2\nu \alpha _{i}}{\alpha _{i}-1}\frac{V_{i}}{R_{i}^{2}}
\label{derv}
\end{equation}
with
\[
\alpha _{i}=\frac{\gamma _{i}+2}{\gamma _{i}}\in \left(
1,\infty\right) .
\]
The pressure $P_{i}$ and the radius $R_{i}$ are related by a function
\[
P_{i}=P_{i}\left( x,R_{i}\right) .
\]
In \cite{SPH02}, it is assumed that
\[
P_{i}\left( x,R_{i}\right) =C\left[ \left( \frac{R_{i}}{R_{0}}\right)
^{\beta }-1\right]
\]
where $C$, $R_{0}$ and $\beta $ are constants. We do not make such an
assumption, only assume that it is a differentiable function and
\[
\frac{\partial P_{i}}{\partial R_{i}}>0
\]
for all $x$ and $R_{i}$. Let
\[
A_{i}=\pi R_{i}^{2},\quad Q_{i}=A_{i}V_{i}
\]
be the cross-section area and the flow rate, respectively. It can
be shown that the system of equations (\ref{massbal}) and
(\ref{derv}) is equivalent to
\begin{equation}
\begin{array}{rcl}
\displaystyle\frac{\partial A_{i}}{\partial t}+\frac{\partial Q_{i}}{%
\partial x} & = & 0, \\[12pt]
\displaystyle\frac{\partial Q_{i}}{\partial t}+\alpha _{i}\frac{\partial }{%
\partial x}\left( \frac{Q_{i}^{2}}{A}\right) +\frac{A_{i}}{\rho }\frac{%
\partial P_{i}}{\partial x} & = & \displaystyle-\frac{4\pi \nu \alpha _{i}}{%
\alpha _{i}-1}\frac{Q_{i}}{A_{i}}.
\end{array}
\label{deaq}
\end{equation}
Also, by rescaling the spatial variable $x$, we may assume that
each vessel is parameterized to $x\in \left( 0,1\right) $.

The system of differential equations are supplemented with the initial
condition
\begin{equation}
P_{i}\left( x,0\right) =P_{i}^{I}\left( x\right) ,\quad Q_{i}\left(
x,0\right) =Q_{i}^{I}\left( x\right)   \label{ic}
\end{equation}
and boundary conditions. Boundary conditions at each end of the
vessel are given according to the type of the end. If it is an
external end of the network, either the pressure
\begin{equation}
P_{i}=P_{i}^{B}\left( t\right)   \label{bcep}
\end{equation}
or the flow rate
\begin{equation}
Q_{i}=Q_{i}^{B}\left( t\right)   \label{bceq}
\end{equation}
is specified. If the end is a \emph{branching junction, }a
junction connecting several vessels, let $j_{1},\ldots ,j_{\nu }$
and $j_{\nu +1},\ldots ,j_{\mu }$ denote the incoming and outgoing
vessels, respectively. One imposes the mass balance condition
\begin{equation}
\sum_{l=1}^{\nu }Q_{j_{l}}\left( 1,t\right) =\sum_{l^{\prime }=\nu +1}^{\mu
}Q_{j_{l^{\prime }}}\left( 0,t\right) ,  \label{bcbq}
\end{equation}
and the momentum balance condition
\begin{equation}
\rho _{j_{l}}\frac{\partial Q_{j_{l}}}{\partial t}=A_{j_{l}}\left(
P_{j_{l}}-P_{junc}\right) ,\quad \rho _{j_{l^{\prime }}}\frac{\partial
Q_{j_{l^{\prime }}}}{\partial t}=A_{j_{l^{\prime }}}\left(
P_{junc}-P_{j_{l^{\prime }}}\right)   \label{bcbp}
\end{equation}
for $l=1\ldots ,\nu $, $l^{\prime }=\nu +1,\ldots ,\mu $, where
$\rho _{i}$ are small positive constants and $P_{junc}$ is the
pressure at the junction. (In \cite{SPH02}, every branching
junction connects only three vessels, it is either a bifurcation
point of one artery into two smaller ones or a joining point of
two veins into a bigger one. Our prescription allows more general
configuration of the network, including the presence of Willis
loops.) If the end is a \emph{transitional junction}, which
connects the vessel to a
network of arterioles, capillaries and venules, we follow the so called\emph{%
\ microcirculation} model proposed in \cite{SPH02,SBL81}.
Generalizing from \cite {SPH02}, arterioles or venules connected
to the vessel $j_{l}$ are represented by a lumped resistive
element $R_{j_{l}}$. The capillary bed is also represented by a
resistive element $R_{C}$. $R_{j_{l}}$'s are connected to $R_{C}$
through capacitive elements $C_{1}$ and $C_{2}$ on the two ends.
\begin{center}
\includegraphics{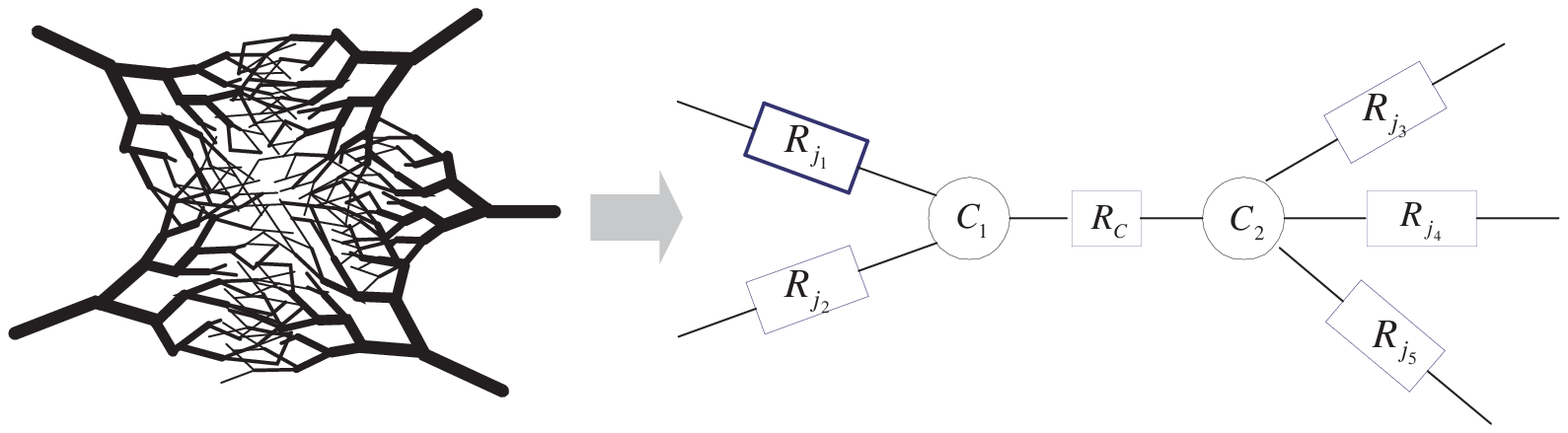}\makeatletter
\makeatletter\def\@captype{figure}\makeatother
\caption{Microcirculation model of the network of arterioles,
capillaries and venules.}
\end{center}
\noindent Let $%
j_{1},\ldots ,j_{\nu }$ be the arteries and let $j_{\nu +1},\ldots
,j_{\mu }$ be the veins that are connected to a
arteriole-capillary-venule network. The boundary conditions for
$P_{j_{l}}$, $Q_{j_{l}}$ are
\begin{equation}
\begin{array}{ll}
R_{j_{l}}Q_{j_{l}}\left( 1,t\right) =P_{j_{l}}\left( 1,t\right)
-P_{C_{1}}\left( t\right) , & \ \text{for }l=1,\ldots ,\nu , \\[8pt]
R_{j_{l^{\prime }}}Q_{j_{l^{\prime }}}\left( 0,t\right) =P_{C_{2}}\left(
0,t\right) -P_{j_{l^{\prime }}}\left( t\right) , & \ \text{for }l^{\prime
}=\nu +1,\ldots ,\mu
\end{array}
\label{bctp}
\end{equation}
and
\begin{equation}
C_{1}\frac{dP_{C_{1}}}{dt}=\sum_{l=1}^{\nu }Q_{j_{l}}\left( 1,t\right)
-Q_{C},\quad C_{2}\frac{dP_{C_{2}}}{dt}=Q_{C}\left( t\right)
-\sum_{l^{\prime }=\nu +1}^{\mu }Q_{j_{l^{\prime }}}\left( 0,t\right)
\label{bctq}
\end{equation}
where $P_{C_{i}}$, $i=1,2$ represent the pressure in the
capacitive elements $C_{1}$, $C_{2}$, and
\begin{equation}
Q_{C}=\frac{P_{C_{1}}-P_{C_{2}}}{R_{C}}  \label{qc}
\end{equation}
represents the flow rate in the resistive element $R_{C}$. (In
\cite {SPH02}, there is only one artery and one vein connected to
the system of arteriole-capillary-venule at the two ends. We do
not rule out the possibility of multiple arteries and veins join
together to such a system.)

The system we study in this paper consists of the equations
\begin{equation}
\begin{array}{ll}
\begin{array}{l}
\displaystyle\frac{\partial P_{i}}{\partial t}+a_{i}\frac{\partial Q_{i}}{%
\partial x}=f_{i}, \\[12pt]
\displaystyle\frac{\partial Q_{i}}{\partial t}+b_{i}\frac{\partial P_{i}}{%
\partial x}+2c_{i}\frac{\partial Q_{i}}{\partial x}=g_{i},
\end{array}
& \quad x\in \left( 0,1\right) ,\ t>0
\end{array}
\label{depq}
\end{equation}
and the initial and boundary conditions given by (\ref{ic})--(\ref{qc}%
). For convenience, we also use the vector form
\begin{equation}
\left( U_{i}\right) _{t}+B_{i}\left( U_{i}\right) _{x}=F_{i}  \label{devec}
\end{equation}
where $U_{i}=\left( P_{i},Q_{i}\right) $, $F_{i}=\left(
f_{i},g_{i}\right) $ and
\[
B_{i}=\left(
\begin{array}{ll}
0 & a_{i} \\
b_{i} & 2c_{i}
\end{array}
\right) .
\]
Eq. (\ref{deaq}) is a special case of this system where
\[
a_{i}=\frac{\partial P_{i}}{\partial A_{i}},\quad b_{i}=\frac{A_{i}}{\rho }-%
\frac{\alpha Q_{i}^{2}}{A_{i}^{2}}\left( \frac{\partial P_{i}}{\partial A_{i}%
}\right) ^{-1},\quad c_{i}=\frac{\alpha Q_{i}}{A_{i}},\ f_{i}=0,\ g_{i}=%
\frac{\alpha Q_{i}^{2}}{A_{i}^{2}}\frac{\partial A_{i}}{\partial x}-\frac{%
4\pi \nu \alpha }{\alpha -1}\frac{Q_{i}}{A_{i}}.
\]
We do not assume any particular form of these functions though,
they are general differentiable functions of $\left(
x,t,P_{i},Q_{i}\right) $. Our basic assumptions are $a_{i}>0$ and
$A_{i}>\varepsilon _{0}$ for some positive constant $\varepsilon
_{0}$. Other assumptions will follow. Apart from the junction
conditions, this system is the same as the one we study in
\cite{RCZC}. Also, the junction conditions in \cite{RCZC} is the
special case of (\ref{bcbq})--( \ref{bcbp}) above with $\rho
_{i}=0$. As in \cite{RCZC}, we use a fixed point principle to
prove the solvability of the
problem. Substituting a pair of functions $\left( p_{i},q_{i}\right) $ for $%
\left( P_{i},Q_{i}\right) $ in the coefficients $a_{i}$, $b_{i}$, $c_{i}$, $%
A_{i}$ and forcing functions $f_{i}$, $g_{i}$, the system becomes linear.
That is, all the functions $a_{i}$, etc. are independent of unknowns. If the
linear system has a unique solution, then, one can establish a mapping from $%
\left( p_{i},q_{i}\right) $ to the linear problem solution $\left(
P_{i},Q_{i}\right) $. If one also shows that this mapping has a
unique fixed point, then the fixed point is necessarily the unique
solution of the quasilinear system. Hence, we shall first give a
condition for the linear system to have a unique solution, then
examine under what conditions the mapping has a unique fixed
point. The first aspect of the problem is investigated in Section
2 and the second in Section 3. We also prove a result on the
continuity of dependence of solutions on the initial, boundary and
forcing functions for linear and quasilinear systems, thus,
completing the analysis of the well-posedness of the problem. In
spite of similarity in parts of the analysis to the one used in
\cite{RCZC}, the more general branching junction condition and the
new transitional junction conditions require more careful
treatments. Hence, there are substantial variations in the
analysis. For completeness and to benefit the reader, we include
all the major arguments in this paper.

\section{The linear system}

\setcounter{equation}{0} \setcounter{theorem}{0}In this section, we analyze (%
\ref{depq}) as a linear system with $a_{i}$, $b_{i}$, $c_{i}$, $f_{i}$, $%
g_{i}$, $A_{i}$ independent of $P_{i}$ and $Q_{i}$. The initial and boundary
conditions are given by (\ref{ic})--(\ref{qc}) except that the junction
condition (\ref{bcbp}) is substituted by the more general condition
\begin{equation}
\rho _{j_{l}}\frac{\partial Q_{j_{l}}}{\partial t}=A_{j_{l}}\left(
P_{j_{l}}-P_{junc}\right) +C_{j_{l}},\quad \rho _{j_{l^{\prime }}}\frac{%
\partial Q_{j_{l^{\prime }}}}{\partial t}=A_{j_{l^{\prime }}}\left(
P_{junc}-P_{j_{l^{\prime }}}\right) +C_{j_{l^{\prime }}}  \label{bcbp'}
\end{equation}
where $C_{i}$ are differentiable functions of $\left( x,t\right)
$. The inclusion of $C_{i}$ is needed in the next section in order
that the result of this section can be extended to the quasilinear
system. We give conditions for the linear system to have a unique
global solution. The conditions are most naturally given in terms
of the eigenvalues of the matrix $B_{i}$, which have the form
\[
\lambda _{i}^{R}=c_{i}+u_{i},\quad \lambda _{i}^{L}=c_{i}-u_{i},
\]
where
\[
u_{i}=\sqrt{c_{i}^{2}+a_{i}b_{i}}.
\]
These eigenvalues are real if
\begin{equation}
c_{i}^{2}+a_{i}b_{i}>0,\ x\in \left( 0,1\right) ,\ t>0,\ i=1,\ldots ,n.
\label{cond2}
\end{equation}
In this case,
\begin{equation}
\lambda _{i}^{R}\left( x,t\right) >0,\text{ }\lambda _{i}^{L}\left(
x,t\right) <\lambda _{i}^{R}\left( x,t\right)   \label{linth1.9}
\end{equation}
and the system is hyperbolic. Under the condition (\ref{cond2}),
we show that the linear system has a unique solution if
\[
\lambda _{i}^{L}\left( 0,t\right) <0,\ \lambda _{i}^{L}\left( 1,t\right)
<0,\ i=1,\ldots ,n
\]
which is equivalent to
\begin{equation}
a_{i}b_{i}>0,\ t\geq 0,\ i=1,\ldots ,n.  \label{cond3}
\end{equation}
at $x=0,1$ only. It needs not hold for $x\in \left( 0,1\right) $.

\begin{theorem}
\smallskip \label{linth1}Assume that the functions $a_{i}$, $b_{i}$, $c_{i}$%
, $f_{i}$, $g_{i}$, $A_{i}$ and $C_{i}$ are independent of $\left(
P_{i},Q_{i}\right) $. Suppose that these functions and the initial and
boundary functions $P_{i}^{I}$, $Q_{i}^{I}$, $P_{i}^{B}$, $Q_{i}^{B}$ all
have bounded first-order derivatives. Suppose also that $a_{i}>0$, $A_{i}>0$
and that the conditions (\ref{cond2}) and (\ref{cond3}) hold. Then, for any $%
T>0$ there is a unique solution in a bounded subset of the space $C\left(
\left[ 0,1\right] \times \left[ 0,T\right] ,\Bbb{R}^{2n}\right) $ to the
linear system (\ref{depq}) with the initial and boundary conditions given by
(\ref{ic})--(\ref{bcbq}), (\ref{bctp})--(\ref{qc}), and (\ref{bcbp'}).
\end{theorem}

\paragraph{Proof.}

\smallskip We first show that the system has a unique solution for $%
0<t<\delta $ for some $\delta >0$. The proof is based on the method of
characteristics and a fixed point principle. For systems defined on only one
branch with boundary conditions of the forms of (\ref{bcep}) or (\ref{bceq}%
), this is a standard approach. In our case, special care is needed to
handle the junction conditions.

Consider the $i$-th branch. From any point $\left( \xi,
\tau\right) $ on the left, right, and lower boundary of the
rectangle $D=:\left[ 0,1\right] \times \left[ 0,T\right] $, we
construct the left-going and right-going
characteristic curves $x=x_{i}^{L}\left( t;\xi ,\tau \right) $ and $%
x=x_{i}^{R}\left( t;\xi ,\tau \right) $ by
\begin{eqnarray*}
\frac{dx_{i}^{L}}{dt} &=&\lambda _{i}^{L}\left( x_{i}^{L},t\right) ,\
x_{i}^{L}\left( \tau \right) =\xi , \\
\frac{dx_{i}^{R}}{dt} &=&\lambda _{i}^{R}\left( x_{i}^{R},t\right) ,\
x_{i}^{R}\left( \tau \right) =\xi ,
\end{eqnarray*}
respectively, where $\lambda _{i}^{L}$ and $\lambda _{i}^{R}$ are the two
eigenvalues of the matrix $B_{i}$. By the uniqueness of solutions to these
differential equations, a left-going (resp. right-going) characteristic
curve cannot intersect with another left-going (resp. right-going)
characteristic curve. Let $X_{i}^{L}$ and $X_{i}^{R}$ be the right-most
left-going and left-most right-going characteristic curves,
\[
x=x_{i}^{L}\left( t;1,0\right) \ \text{and }x=x_{i}^{R}\left( t;0,0\right)
\]
starting from the lower boundary of $D$, respectively. It can be shown from (%
\ref{linth1.9}) that the two curves can have at most one intersection. Let $%
t_{i}$ be the value of $t$ at the intersection. If the two curves do not
intersect in $D$, we simply define $t_{i}=T$. By condition (\ref{cond3}), $%
X_{i}^{L}$ cannot reach the right vertical line $x=1$ at any $t>0$, and by $%
\lambda _{i}^{R}>0$, $X_{i}^{R}$ cannot reach the vertical line $x=0$ at any
$t>0$. Thus, the rectangle $D_{i}=:\left[ 0,1\right] \times \left[
0,t_{i}\right] $ can be divided into three parts
\[
D_{i}=D_{i}^{L}\cup D_{i}^{C}\cup D_{i}^{R},
\]
where $D_{i}^{L}$ is between the vertical line $x=0$ and the characteristic
curve $X_{i}^{R}$, $D_{i}^{C}$ is between the two characteristic curves, and
$D_{i}^{R}$ is between $X_{i}^{L}$ and $x=1$.

\begin{center}
\includegraphics{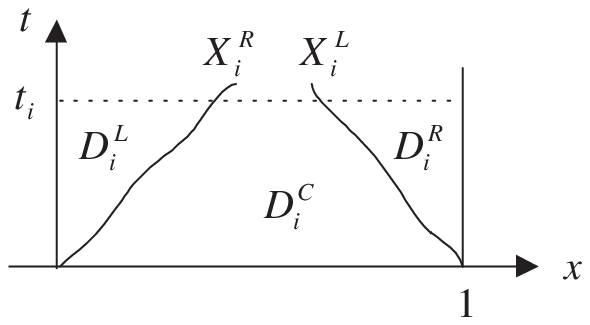}\makeatletter
\makeatletter\def\@captype{figure}\makeatother \caption{Three
parts of $D_{i}$}
\end{center}

\noindent We show that there is a $\delta _{i}\leq t_{i}$ such that the
solution $\left( P_{i},Q_{i}\right) $ for the $i$-th branch exists in the
restriction of $D_{i}$ to the strip $\left\{ 0\leq t\leq \delta _{i}\right\}
$.

First, observe that the initial conditions alone determine the solution
completely in the central region $D_{i}^{C}$. This follows from the theory
of first-order linear hyperbolic systems and the fact that from any point $%
\left( x,t\right) \in D_{i}^{C}$, the two characteristic curves, followed
backwards, must land on the horizontal line $t=0$. (The latter is a
consequence of (\ref{linth1.9}).) To extend the solution to other parts of $%
D_{i}$, we make a change of unknowns and derive a set of integral equations.
Note that $l_{i}^{R}=:\left( -\lambda _{i}^{L},a_{i}\right) $ and $%
l_{i}^{L}=:\left( -\lambda _{i}^{R},a_{i}\right) $ are the left eigenvectors
of $B_{i}$ corresponding to $\lambda _{i}^{R}$ and $\lambda _{i}^{L}$,
respectively. Introduce new unknowns
\begin{equation}
r_{i}=l_{i}^{R}U_{i}\equiv -\lambda _{i}^{L}P_{i}+a_{i}Q_{i},\quad
s_{i}=l_{i}^{L}U_{i}\equiv -\lambda _{i}^{R}P_{i}+a_{i}Q_{i}.
\label{linth1.6}
\end{equation}
The system (\ref{depq}) can be written in terms of $r_{i}$ and $%
s_{i}$ by multiplying the left eigenvectors to (\ref{devec}) and
substituting in
\begin{equation}
P_{i}=\frac{1}{2u_{i}}\left( r_{i}-s_{i}\right) ,\quad Q_{i}=\frac{1}{%
2u_{i}a_{i}}\left( \lambda _{i}^{R}r_{i}-\lambda _{i}^{L}s_{i}\right) .
\label{linth1.1}
\end{equation}
This results in the equations
\begin{equation}
\partial _{i}^{R}r_{i}=F_{i}^{R}\left( x,t,r_{i},s_{i}\right) ,\ \partial
_{i}^{L}s_{i}=F_{i}^{L}\left( x,t,r_{i},s_{i}\right) ,  \label{linth1.2}
\end{equation}
where
\begin{equation}
\partial _{i}^{R}=\frac{\partial }{\partial t}+\lambda _{i}^{R}\frac{%
\partial }{\partial x},\quad \partial _{i}^{L}=\frac{\partial }{\partial t}%
+\lambda _{i}^{L}\frac{\partial }{\partial x},  \label{linth1.5}
\end{equation}
and
\begin{equation}
F_{i}^{R}\left( x,t,r_{i},s_{i}\right) =l_{i}^{R}F_{i}+\left( \partial
_{i}^{R}l_{i}^{R}\right) U_{i},\quad F_{i}^{L}\left( x,t,r_{i},s_{i}\right)
=l_{i}^{L}F_{i}+\left( \partial _{i}^{L}l_{i}^{L}\right) U_{i}.
\label{linth1.3}
\end{equation}
(A differential operator acting on a vector means that it acts on
each component of the vector.) Let $\left( x,t\right) \in D_{i}$.
We integrate the first equation of (\ref{linth1.2}) along the
right-going characteristic curve $x^{R}\left( t;\xi ,\tau \right)
$ which passes through
$\left( x,t\right) $ and reaches the left or lower boundary of $D_{i}$ at $%
\left( \xi ,\tau \right) $. It can be shown that for $\left( x,t\right) \in
D_{i}^{C}\cup D_{i}^{R}$, $\tau =0$, and for $\left( x,t\right) \in D_{i}^{L}
$, $\xi =0$. In the former case, we obtain
\begin{equation}
r_{i}\left( x,t\right) =r_{i}^{I}\left( \xi \right)
+\int_{0}^{t}F_{i}^{R}\left( x_{i}^{R}\left( t^{\prime };\xi ,0\right)
,t^{\prime },r_{i},s_{i}\right) dt^{\prime }  \label{inteq1}
\end{equation}
In the latter case, we have
\begin{equation}
r_{i}\left( x,t\right) =r_{i}\left( 0,\tau \right) +\int_{\tau
}^{t}F_{i}^{R}\left( x_{i}^{R}\left( t^{\prime };0,\tau \right) ,t^{\prime
},r_{i},s_{i}\right) dt^{\prime }.  \label{inteq2}
\end{equation}
Similarly, by integrating the second equation of (\ref{linth1.2})
along the left-going characteristic curve $x_{i}^{L}\left( t;\xi
,\tau \right) $ that passes through both $\left( x,t\right) $ and
$\left( \xi ,\tau \right) $ (which is on either the right or lower
boundary of $D_{i}$), the equations are
\begin{equation}
s_{i}\left( x,t\right) =s_{i}^{I}\left( \xi \right)
+\int_{0}^{t}F_{i}^{L}\left( x_{i}^{L}\left( t^{\prime };\xi ,0\right)
,t^{\prime },r_{i},s_{i}\right) dt^{\prime }  \label{inteq3}
\end{equation}
if $\left( x,t\right) \in D_{i}^{L}\cup D_{i}^{C}$ and
\begin{equation}
s_{i}\left( x,t\right) =s_{i}\left( 1,\tau \right) +\int_{\tau
}^{t}F_{i}^{L}\left( x_{i}^{L}\left( t^{\prime };1,\tau \right) ,t^{\prime
},r_{i},s_{i}\right) dt^{\prime }  \label{inteq4}
\end{equation}
if $\left( x,t\right) \in D_{i}^{R}$. These are the integral
equations we need.

For any $\delta _{i}\leq t_{i}$ we use $D_{i,\delta _{i}}^{L}$, $D_{i,\delta
_{i}}^{C}$ and $D_{i,\delta _{i}}^{R}$ to denote the restrictions of $%
D_{i}^{L}$, $D_{i}^{C}$ and $D_{i}^{R}$ to the strip $\left\{ 0\leq t\leq
\delta _{i}\right\} $, respectively. First, consider the case where the end
of the branch is an external end. We discuss the case of a left end only,
the case of a right end can be treated similarly. If the boundary condition
is given by (\ref{bcep}), we define $\hat{s}_{i}=s_{i}/\varepsilon $ where $%
\varepsilon <1$ is any constant. Using the first equation of (\ref{linth1.1}%
) in the integral equations (\ref{inteq2}) and (\ref{inteq3}),
\begin{equation}
\left(
\begin{array}{l}
r_{i}\left( x,t\right)  \\
\hat{s}_{i}\left( x,t\right)
\end{array}
\right) =\left(
\begin{array}{l}
\displaystyle2u_{i}\left( 0,\tau \right) P_{i}^{B}\left( \tau \right)
+\varepsilon \hat{s}_{i}\left( 0,\tau \right) +\int_{\tau
}^{t}F_{i}^{R}\left( x_{i}^{R}\left( t^{\prime };0,\tau \right) ,t^{\prime
},r_{i},\varepsilon \hat{s}_{i}\right) dt^{\prime } \\
\displaystyle\frac{1}{\varepsilon }s_{i}^{I}\left( \xi \right) +\frac{1}{%
\varepsilon }\int_{0}^{t}F_{i}^{L}\left( x_{i}^{L}\left( t^{\prime };\xi
,0\right) ,t^{\prime },r_{i},\varepsilon \hat{s}_{i}\right) dt^{\prime }
\end{array}
\right) .  \label{inteq5}
\end{equation}
This is a fixed point equation for $\left(
r_{i},\hat{s}_{i}\right)
$ if we define the right hand side as a mapping of an operator $K$ on $%
\left( r_{i},\hat{s}_{i}\right) $ in a bounded subset of $C\left(
D_{i,\delta _{i}}^{L}\cup D_{i,\delta _{i}}^{C},\Bbb{R}^{2}\right) $. In a
standard approach, it can be shown that $K$ is a contraction mapping if $%
\delta _{i}$ is sufficiently small. Hence, the fixed point exists and is
unique, and the solution $\left( r_{i},s_{i}\right) $ is uniquely extended
to $D_{i,\delta _{i}}^{L}\cup D_{i,\delta _{i}}^{C}$. If the boundary
condition is given by (\ref{bceq}), we define $\hat{s}_{i}=s_{i}/\varepsilon
$, where $\varepsilon >0$ is so small such that
\[
\varepsilon \left| \frac{\lambda _{i}^{L}\left( 0,\tau \right) }{\lambda
_{i}^{R}\left( 0,\tau \right) }\right| <1,\quad \tau \in \left(
0,t_{i}\right) .
\]
The fixed point equation is then
\begin{equation}
\left(
\begin{array}{l}
r_{i}\left( x,t\right)  \\
\hat{s}_{i}\left( x,t\right)
\end{array}
\right) =\left(
\begin{array}{l}
\displaystyle\frac{2a_{i}u_{i}\left( 0,\tau \right) }{\lambda _{i}^{R}\left(
0,\tau \right) }Q_{i}^{B}\left( \tau \right) +\frac{\lambda _{i}^{L}\left(
0,\tau \right) }{\lambda _{i}^{R}\left( 0,\tau \right) }\varepsilon \hat{s}%
_{i}\left( 0,\tau \right) +\int_{\tau }^{t}F_{i}^{R}\left( x_{i}^{R}\left(
t^{\prime };0,\tau \right) ,t^{\prime },r_{i},\varepsilon \hat{s}_{i}\right)
dt^{\prime } \\
\displaystyle\frac{1}{\varepsilon }s_{i}^{I}\left( \xi \right) +\frac{1}{%
\varepsilon }\int_{0}^{t}F_{i}^{L}\left( x_{i}^{L}\left( t^{\prime };\xi
,0\right) ,t^{\prime },r_{i},\varepsilon \hat{s}_{i}\right) dt^{\prime }
\end{array}
\right) .  \label{inteq5'}
\end{equation}
By a similar argument, the solution can again be uniquely
extended.

We next extend the solution to either $D_{i,\delta _{i}}^{L}$ or $%
D_{i,\delta _{i}}^{R}$ if the end is a branching junction. In this case, we
shall extend the solution on all the branches that are connected to the same
junction simultaneously. Let $j_{1},\ldots ,j_{\nu }$ be the incoming and $%
j_{\nu +1},\ldots ,j_{\mu }$ the outgoing branches to the junction.
Equations (\ref{bcbq}), (\ref{bcbp'}) and (\ref{linth1.1}) give rise to a $%
2\mu \times \mu $ homogenous system of linear (ordinary) differential
equations for $r_{i}\left( 1,t\right) $, $s_{i}\left( 1,t\right) $, $%
i=j_{1},\ldots ,j_{\nu }$ and $r_{i}\left( 0,t\right) $, $s_{i}\left(
0,t\right) $, $i=j_{\nu +1},\ldots ,j_{\mu }$:
\begin{equation}
\begin{array}{l}
\displaystyle\frac{\rho _{j_{1}}}{A_{j_{1}}}\frac{d}{dt}Q_{j_{1}}\left(
1,t\right) -\frac{C_{j_{1}}}{A_{j_{1}}}-P_{j_{1}}\left( 1,t\right) =\frac{%
\rho _{i}}{A_{i}}\frac{d}{dt}Q_{i}\left( 1,t\right) -\frac{C_{i}}{A_{i}}%
-P_{i}\left( 1,t\right) ,\quad i=j_{2},\ldots ,j_{\nu }, \\[12pt]
\displaystyle\frac{\rho _{j_{1}}}{A_{j_{1}}}\frac{d}{dt}Q_{j_{1}}\left(
1,t\right) -\frac{C_{j_{1}}}{A_{j_{1}}}-P_{j_{1}}\left( 1,t\right) =-\frac{%
\rho _{i}}{A_{i}}\frac{d}{dt}Q_{i}\left( 0,t\right) -\frac{C_{i}}{A_{i}}%
-P_{i}\left( 0,t\right) ,\quad i=j_{\nu +1},\ldots ,j_{\mu }, \\[12pt]
\displaystyle\sum_{l=1}^{\nu }Q_{j_{l}}\left( 1,t\right) -\sum_{l^{\prime
}=\nu +1}^{\mu }Q_{j_{l^{\prime }}}\left( 0,\tau \right) =0.
\end{array}
\label{linth1.10}
\end{equation}
Differentiate the last equation with respect to $t$ and regard $%
s_{j_{1}}\left( 1,t\right)$,..., $s_{j_{\nu }}\left( 1,t\right) $, $%
r_{j_{\nu +1}}\left( 0,t\right)$,..., $r_{j_{\mu }}\left(
0,t\right) $ as unknowns. The derivatives of unknowns can be
solved from (\ref{linth1.10}) because the coefficient matrix of
$ds_{i}/dt$ and $dr_{i}/dt$ in (\ref{linth1.10}),
\[
\left(
\begin{array}{cccc}
-\frac{\rho _{j_{1}}\lambda _{j_{1}}^{L}\left( 1,t\right) }{%
2u_{j_{1}}a_{j_{1}}A_{j_{1}}\left( 1,t\right) } & \frac{\rho _{j_{2}}\lambda
_{j_{2}}^{L}\left( 1,t\right) }{2u_{j_{2}}a_{j_{2}}A_{j_{2}}\left(
1,t\right) } & \cdots  & 0 \\
\vdots  & \vdots  & \ddots  & \vdots  \\
-\frac{\rho _{j_{1}}\lambda _{j_{1}}^{L}\left( 1,t\right) }{%
2u_{j_{1}}a_{j_{1}}A_{j_{1}}\left( 1,t\right) } & 0 & \cdots  & \frac{\rho
_{j_{\mu }}\lambda _{j_{\mu }}^{R}\left( 0,t\right) }{2u_{j_{\mu }}a_{j_{\mu
}}A_{j_{\mu }}\left( 0,t\right) } \\
-\frac{\lambda _{j_{1}}^{L}\left( 1,t\right) }{2u_{j_{1}}a_{j_{1}}\left(
1,t\right) } & -\frac{\lambda _{j_{2}}^{L}\left( 1,t\right) }{%
2u_{j_{2}}a_{j_{2}}\left( 1,t\right) } & \cdots  & -\frac{\lambda _{j_{\mu
}}^{R}\left( 0,t\right) }{2u_{j_{\mu }}a_{j_{\mu }}\left( 0,t\right) }
\end{array}
\right)
\]
has the determinant
\[
\left( \frac{-1}{2}\right) ^{\mu }\prod_{l=1}^{\nu }\frac{\rho
_{j_{l}}\lambda _{j_{l}}^{L}\left( 1,t\right) }{u_{j_{l}}a_{j_{l}}A_{j_{l}}%
\left( 1,t\right) }\prod_{l^{\prime }=\nu +1}^{\mu }\frac{\rho
_{j_{l^{\prime }}}\lambda _{j_{l^{\prime }}}^{R}\left( 0,t\right) }{%
u_{j_{l^{\prime }}}a_{j_{l^{\prime }}}A_{j_{l^{\prime }}}\left( 0,t\right) }%
\sum_{l=1}^{\mu }\frac{A_{j_{l}}}{\rho _{j_{l}}}.
\]
Since $\lambda _{i}^{L}<0<\lambda _{i}^{R}$ at the junction, the determinant
is not zero. Thus, the derivatives of the unknowns, $s_{j_{l}}\left(
1,t\right) $ and $r_{j_{l^{\prime }}}\left( 0,t\right) $, are each a linear
combination of the functions $r_{j_{l}}\left( 1,t\right) $, $s_{j_{l}}\left(
1,t\right) $, $r_{j_{l^{\prime }}}\left( 0,t\right) $, $s_{j_{l^{\prime
}}}\left( 0,t\right) $ together with the derivatives of $r_{j_{l}}\left(
1,t\right) $ and $s_{j_{l^{\prime }}}\left( 0,t\right) $, $l=1,\ldots ,\nu $%
, $l^{\prime }=\nu +1,\ldots ,\mu $. Integrating and using the initial
condition determined by (\ref{ic}) and (\ref{linth1.6}), we can write
\begin{equation}
s_{i}\left( 1,\tau \right) =s_{i}\left( 1,0\right) +\sum_{l=1}^{\nu
}m_{j_{l}}^{i}\left( \tau \right) r_{j_{l}}\left( 1,\tau \right)
+\sum_{l^{\prime }=\nu +1}^{\mu }m_{j_{l^{\prime }}}^{i}\left( \tau \right)
s_{j_{l^{\prime }}}\left( 0,\tau \right) +\int_{0}^{\tau }H_{i}dt^{\prime },
\label{linth1.11}
\end{equation}
for $i=j_{1},\ldots ,j_{\nu }$ and
\begin{equation}
r_{i}\left( 0,\tau \right) =r_{i}\left( 0,0\right) +\sum_{l=1}^{\nu
}n_{j_{l}}^{i}\left( \tau \right) r_{j_{l}}\left( 1,\tau \right)
+\sum_{l^{\prime }=\nu +1}^{\mu }n_{j_{l^{\prime }}}^{i}\left( \tau \right)
s_{j_{l^{\prime }}}\left( 0,\tau \right) +\int_{0}^{\tau }H_{i}dt^{\prime },
\label{linth1.12}
\end{equation}
for $i=j_{\nu +1},\ldots ,j_{\mu }$, where $m_{j}^{i}$, $n_{j}^{i}$ are
continuous functions and $H_{i}$ are linear combinations of $r_{j_{l}}\left(
1,t\right) $, $s_{j_{l}}\left( 1,t\right) $, $r_{j_{l^{\prime }}}\left(
0,t\right) $, $s_{j_{l^{\prime }}}\left( 0,t\right) $ and $C_{i}\left(
t\right) $ with coefficients depending only on $t$. Choose an $\varepsilon >0
$ such that
\[
\varepsilon \max \left\{ \sum_{l=1}^{\mu }\left| m_{j_{l}}^{i}\left( \tau
\right) \right| ,\sum_{l=1}^{\mu }\left| n_{j_{l}}^{i}\left( \tau \right)
\right| \right\} <1,\quad i=j_{1},\ldots ,j_{\mu },\ \tau \in \left[
0,t_{i}\right]
\]
and introduce
\[
\hat{r}_{j_{l}}=\frac{r_{j_{l}}}{\varepsilon },\quad \hat{s}_{j_{l^{\prime
}}}=\frac{s_{j_{l^{\prime }}}}{\varepsilon },\quad l=1,\ldots ,\nu ,\
l^{\prime }=\nu +1,\ldots ,\mu .
\]
Then, from (\ref{inteq1})--(\ref{inteq4}), the integral equations for the $%
2\mu $ unknowns $\hat{r}_{j_{l}}$, $s_{j_{l}}$, $r_{j_{l^{\prime }}}$, $\hat{%
s}_{j_{l^{\prime }}}$, $l=1,\ldots ,\nu $, $l^{\prime }=\nu +1,\ldots ,\mu $
constitute a fixed point equation, $w=Kw$, where
\begin{equation}
w=\left( \hat{r}_{j_{1}},\ldots ,\hat{r}_{j_{\nu }},s_{j_{1}},\ldots
,s_{j_{\nu }},r_{j_{\nu +1}},\ldots ,r_{j_{\mu }},\hat{s}_{j_{\nu
+1}},\ldots ,\hat{s}_{j_{\mu }}\right)   \label{linth1.7}
\end{equation}
and
\begin{equation}
\begin{array}{l}
Kw=\left( \frac{1}{\varepsilon }r_{j_{1}}^{I}\left( \xi _{j_{1}}\right) +%
\frac{1}{\varepsilon }\int_{0}^{t}F_{j_{1}}^{R}dt^{\prime },\ldots ,\right.
\\
\quad s_{j_{1}}\left( 1,0\right) +\varepsilon \left( \sum_{k=1}^{\nu
}m_{j_{k}}^{1}\hat{r}_{j_{k}}\left( 1,\tau \right) +\sum_{k^{\prime }=\nu
+1}^{\mu }m_{j_{k^{\prime }}}^{1}\hat{s}_{j_{k^{\prime }}}\left( 0,\tau
\right) \right) +\int_{0}^{\tau }H_{j_{1}}dt^{\prime }+\int_{\tau
}^{t}F_{j_{1}}^{L}dt^{\prime },\ldots , \\
\quad r_{j_{\nu +1}}\left( 0,0\right) +\varepsilon \left( \sum_{k=1}^{\nu
}n_{j_{k}}^{1}\hat{r}_{j_{k}}\left( 1,\tau \right) +\sum_{k^{\prime }=\nu
+1}^{\mu }n_{j_{k^{\prime }}}^{1}\hat{s}_{j_{k^{\prime }}}\left( 1,\tau
\right) \right) +\int_{0}^{\tau }H_{j_{\nu +1}}dt^{\prime }+\int_{\tau
}^{t}F_{j_{\nu +1}}^{R}dt^{\prime },\ldots , \\
\quad \left. \frac{1}{\varepsilon }s_{j_{\nu +1}}^{I}\left( \xi _{j_{\nu
+1}}\right) +\frac{1}{\varepsilon }\int_{0}^{t}F_{j_{\nu +1}}^{L}dt^{\prime
},\ldots \right)
\end{array}
\label{linth1.8}
\end{equation}
in which
\[
F_{j_{l}}^{R}=F_{j_{l}}^{R}\left( x_{j_{l}}^{R},t^{\prime },\varepsilon \hat{%
r}_{j_{l}},s_{j_{l}}\right) ,\quad F_{j_{l}}^{L}=F_{j_{l}}^{L}\left(
x_{j_{l}}^{L},t^{\prime },\varepsilon \hat{r}_{j_{l}},s_{j_{l}}\right)
\]
for $l=1,\ldots ,\nu $, and
\[
F_{j_{l^{\prime }}}^{R}=F_{j_{l^{\prime }}}^{R}\left( x_{j_{l^{\prime
}}}^{R},t^{\prime },r_{j_{l^{\prime }}},\varepsilon \hat{s}_{j_{l^{\prime
}}}\right) ,\quad F_{j_{l^{\prime }}}^{L}=F_{j_{l^{\prime }}}^{L}\left(
x_{j_{l^{\prime }}}^{L},t^{\prime },r_{j_{l^{\prime }}},\varepsilon \hat{s}%
_{j_{l^{\prime }}}\right)
\]
for $l^{\prime }=\nu +1,\ldots ,\mu $. It can be shown by a standard
argument that $K$ is a contraction mapping in the space
\[
X_{j}=:\prod_{l=1}^{\nu }C\left( D_{j_{l},\delta _{j}}^{C}\cup
D_{j_{l},\delta _{j}}^{R},\Bbb{R}^{2}\right) \times \prod_{l=\nu +1}^{\mu
}C\left( D_{j_{l},\delta _{j}}^{L}\cup D_{j_{l},\delta _{j}}^{L},\Bbb{R}%
^{2}\right)
\]
if $\delta _{j}$ is sufficiently small. Hence, it has a unique fixed point
in $X_{j}$. This extends the solution $\left( r_{i},s_{i}\right) $ for the
neighboring branches of the junction.

It remains to extend the solution to a region adjacent to a
transitional junction. Similar to the case of a branching
junction, we simultaneousely treat all the branches that are
connected to the same transitional junction. Let $j_{1},\ldots
,j_{\nu }$ be the arteries and $j_{\nu +1},\ldots ,j_{\mu }
$ be the veins. The condition connecting the vessels are given by (\ref{bctp}%
), (\ref{bctq}) and (\ref{qc}). Differentiate the equations in (\ref{bctp})
with respect to $t$, the resulting equations together with (\ref{bctq}) is a
linear system of the derivatives of the functions $r_{j_{l}}\left(
1,t\right) $, $s_{j_{l}}\left( 1,t\right) $, $r_{j_{l^{\prime }}}\left(
0,t\right) $, $s_{j_{l^{\prime }}}\left( 0,t\right) $, $l=1,\ldots ,\nu $, $%
l^{\prime }=\nu +1,\ldots ,\mu $, and $P_{C_{1}}\left( t\right) $, $%
P_{C_{2}}\left( t\right) $. The coefficient matrix of $ds_{j_{l}}\left(
1,t\right) /dt$, $dr_{j_{l^{\prime }}}\left( 0,t\right) /dt$, for $%
l=1,\ldots ,\nu $, $l^{\prime }=\nu +1,\ldots ,\mu $ and $P_{C_{1}}^{\prime }
$, $P_{C_{2}}^{\prime }$ is
\[
\left[
\begin{array}{ccc}
D_{1} & 0 & B_{1} \\
0 & D_{2} & B_{2} \\
0 & 0 & I_{2}
\end{array}
\right]
\]
where
\begin{eqnarray*}
D_{1} &=&\text{diag }\left( -\frac{R_{j_{1}}\lambda _{j_{1}}^{L}}{%
2u_{j_{1}}a_{j_{1}}}+\frac{1}{2u_{j_{1}}},\ldots ,-\frac{R_{j_{\nu }}\lambda
_{j_{\nu }}^{L}}{2u_{j_{\nu }}a_{j_{\nu }}}+\frac{1}{2u_{j_{\nu }}}\right) ,
\\
D_{2} &=&\text{diag }\left( \frac{R_{j_{\nu +1}}\lambda _{j_{\nu +1}}^{R}}{%
2u_{j_{\nu +1}}a_{j_{\nu +1}}}+\frac{1}{2u_{j_{\nu +1}}},\ldots ,\frac{%
R_{j_{\mu }}\lambda _{j_{\mu }}^{R}}{2u_{j_{\mu }}a_{j_{\mu }}}+\frac{1}{%
2u_{j_{\mu }}}\right) ,
\end{eqnarray*}
$I_{2}$ is the $2\times 2$ identity matrix, and $B_{1}$, $B_{2}$ are some
constant matrices. Since all the elements of the diagonal matrices $D_{1}$
and $D_{2}$ are positive, the system can be uniquely solved for these
derivatives. Thus, each of $ds_{j_{l}}\left( 1,t\right) /dt$, $%
dr_{j_{l^{\prime }}}\left( 0,t\right) /dt$ for $l=1,\ldots ,\nu $, $%
l^{\prime }=\nu +1,\ldots ,\mu $ is a linear combination of $r_{j_{l}}\left(
1,t\right) $, $s_{j_{l}}\left( 1,t\right) $, $r_{j_{l^{\prime }}}\left(
0,t\right) $, $s_{j_{l^{\prime }}}\left( 0,t\right) $, and $dr_{j_{l}}\left(
1,t\right) /dt$, $ds_{j_{l^{\prime }}}\left( 0,t\right) /dt$, $l=1,\ldots
,\nu $, $l^{\prime }=\nu +1,\ldots ,\mu $ as well as $P_{C_{1}}$ and $%
P_{C_{2}}$. We can also eliminate $P_{C_{1}}$ and $P_{C_{2}}$ by (\ref{bctp}%
). Integrating the resulting equations, we obtain equations (\ref{linth1.11}%
)--(\ref{linth1.12}) for some functions $m_{j}^{i}$, $n_{j}^{i}$, and $H_{i}$%
. The remaining of the previous paragraph can then be used here to give the
extension of the solution to the left or right regions for the branches.

Finally, if we let $\delta $ be the minimum of all $\delta _{i}$ occurring
above, we see that $\delta >0$ and the solution exists and is unique in $%
\left( x,t\right) \in D_{\delta }=:\left[ 0,1\right] \times \left[ 0,\delta
\right] $. Observe that $\delta $ depends only on the bounds of the system
functions $a_{i}$, etc., the initial and boundary functions $P_{i}^{I}$,
etc., and their first-order derivatives in $D=\left[ 0,1\right] \times
\left[ 0,T\right] $. Hence, it is independent of $t$, and we can extend the
solution successively in the time intervals $\left[ 0,\delta \right] $, $%
\left[ \delta ,2\delta \right] $, etc. In this way, the solution is obtained
in $D$ in finitely many steps. \endproof

\smallskip We next derive an estimate of the deviation of solution in term
of the deviations of the initial, boundary and forcing functions. This
estimate is needed in the next section. For any vector function $v=\left(
v_{1},\ldots ,v_{k}\right) $ defined in $C\left( X;\Bbb{R}^{k}\right) $, we
use $\left| v\right| _{X}$ to denote the norm $\max_{i}\left\{ \left|
v_{i}\right| _{C\left( X\right) }\right\} $, where $X$ represents a closed
subset of either $\Bbb{R}$ or $\Bbb{R}^{2}$.

\begin{lemma}
\label{conth1}Let $U=\left( P,Q\right) $ and $\tilde{U}=\left( \tilde{P},%
\tilde{Q}\right) $ be two solutions of the linear problem (\ref{devec}) with
different initial, boundary, and forcing functions. Suppose the conditions
of Theorem \ref{linth1} hold for both solutions. Suppose also that there is
a positive lower bound for all $A_{i}$. Then, there exists a constant $M>0$,
independent of initial, boundary and forcing functions, such that
\begin{equation}
\begin{array}{r}
\left| U-\tilde{U}\right| _{C\left( D_{\delta }\right) }\leq M\left( \left|
P^{I}-\tilde{P}^{I}\right| _{C\left[ 0,1\right] }+\left| Q^{I}-\tilde{Q}%
^{I}\right| _{C\left[ 0,1\right] }+\left| P^{B}-\tilde{P}^{B}\right|
_{C\left[ 0,\delta \right] }+\left| Q^{B}-\tilde{Q}^{B}\right| _{C\left[
0,\delta \right] }\right.  \\
+\left| C-\tilde{C}\right| _{C\left[ 0,\delta \right] }+\left. \delta \left|
f-\tilde{f}\right| _{C\left( D_{\delta }\right) }+\delta \left| g-\tilde{g}%
\right| _{C\left( D_{\delta }\right) }\right) .
\end{array}
\label{conth1.1}
\end{equation}
\end{lemma}

\paragraph{\noindent Proof.}

We need only prove (\ref{conth1.1}) for a $\delta \leq \min_{i}\left\{
\delta _{i}\right\} $, where $\delta _{i}$ represents the constants
occurring in the proof of Theorem \ref{linth1}. This is because for larger $%
\delta $, we can divide the interval $\left[ 0,\delta \right] $ into
subintervals, each has a length less than $\min_{i}\left\{ \delta
_{i}\right\} $, and apply (\ref{conth1.1}) in each subinterval. We can then
take the maximum on each side of the inequalities to derive the inequality
of in $\left[ 0,\delta \right] $. In the sequel, $D_{\delta }^{C}$, $%
D_{\delta }^{L}$ and $D_{\delta }^{R}$ are the restrictions of $D_{i}^{C}$, $%
D_{i}^{L}$ and $D_{i}^{R}$ to the strip $\left\{ 0\leq t\leq \delta \right\}
$, respectively.

By linearity, $U-\tilde{U}$ is the solution of the system with the initial,
boundary and forcing functions $P_{i}^{I}-\tilde{P}_{i}^{I}$, $Q_{i}^{I}-%
\tilde{Q}_{i}^{I}$, $P_{i}^{B}-\tilde{P}_{i}^{B}$, $Q_{i}^{B}-\tilde{Q}%
_{i}^{B}$, $f_{i}-\tilde{f}_{i}$, $g_{i}-\tilde{g}_{i}$ and $C_{i}-\tilde{C}%
_{i}$. Let $r_{i}$, $\hat{r}_{i}$, $s_{i}$, $\hat{s}_{i}$ be defined as in
the proof of Theorem \ref{linth1}, corresponding to $U-\tilde{U}$. We show
that these quantities have upper bounds in the form of the right hand side
of (\ref{conth1.1}) in $D_{\delta }^{C}$, $D_{\delta }^{L}$ and $D_{\delta
}^{R}$.

In $D_{\delta }^{C}$, (\ref{inteq1}) and (\ref{inteq3}) hold. Notice that
the functions $F_{i}^{R}$ and $F_{i}^{L}$ are linear in $r_{i}$, and $s_{i}$%
. Hence, there exists a constant $M$ (we will use $M$ generically for any
constant bounds that are independent of solutions) such that
\[
R_{i}^{C}\left( t\right) +S_{i}^{C}\left( t\right) \leq \left|
r_{i}^{I}\right| _{C\left[ 0,1\right] }+\left| s_{i}^{I}\right| _{C\left[
0,1\right] }+M\int_{0}^{t}\left( R_{i}^{C}\left( t^{\prime }\right)
+S_{i}^{C}\left( t^{\prime }\right) +T_{i}^{C}\left( t^{\prime }\right)
\right) dt^{\prime },
\]
where
\begin{equation}
R_{i}^{C}\left( t\right) =\sup_{\left\{ x:\left( x,t\right) \in D_{\delta
}^{C}\right\} }\left| r_{i}\left( x,t\right) \right| ,\quad S_{i}^{C}\left(
t\right) =\sup_{\left\{ x:\left( x,t\right) \in D_{\delta }^{C}\right\}
}\left| s_{i}\left( x,t\right) \right| ,  \label{conth1.3}
\end{equation}
and
\begin{equation}
T_{i}^{C}\left( t\right) =\sup_{\left\{ x:\left( x,t\right) \in D_{\delta
}^{C}\right\} }\left( \left| f_{i}\left( x,t\right) -\tilde{f}_{i}\left(
x,t\right) \right| +\left| g_{i}\left( x,t\right) -\tilde{g}_{i}\left(
x,t\right) \right| \right) .  \label{conth1.4}
\end{equation}
Hence, by Gronwall's inequality (see, e.g. \cite[p.327]{Mco95}),
\[
R_{i}^{C}\left( t\right) +S_{i}^{C}\left( t\right) \leq M\left( \left|
r_{i}^{I}\right| _{C\left[ 0,1\right] }+\left| s_{i}^{I}\right| _{C\left[
0,1\right] }+\delta \sup_{t\in \left( 0,\delta \right) }T_{i}^{C}\left(
t\right) \right)
\]
for $t\in \left[ 0,\delta \right] $. This proves that $R_{i}^{C}$ and $%
S_{i}^{C}$ have upper bounds in the form of the right side of (\ref{conth1.1}%
).

We next consider the left or right regions if it is adjacent to an external
end. Since the both cases are similar, we will only treat the case where the
left end is an external. The integral equations to be used are (\ref{inteq5}%
) or (\ref{inteq5'}) according to the type of the boundary condition. The
resulting inequality has the form
\[
R_{i}^{L}\left( t\right) +\hat{S}_{i}^{L}\left( t\right) \leq \sigma \hat{S}%
_{i}^{L}\left( t\right) +M\left( \left| s_{i}^{I}\right| _{C\left[
0,1\right] }+\left| \xi _{i}^{B}\right| _{C\left[ 0,\delta \right]
}+\int_{0}^{t}\left( R_{i}^{L}\left( \tau \right) +\hat{S}_{i}^{L}\left(
\tau \right) +T_{i}^{L}\left( \tau \right) \right) d\tau \right)
\]
where $\xi _{i}^{B}$ is either $P_{i}^{B}$ or $Q_{i}^{B}$ depending on the
boundary condition, and $R_{i}^{L}$, $\hat{S}_{i}^{L}$ and $T_{i}^{L}$ are
defined in the same way as in (\ref{conth1.3})--(\ref{conth1.4}), with $%
D_{\delta }^{C}$ substituted by $D_{\delta }^{L}\cup D_{\delta }^{C}$, and $%
\sigma >0$ is a positive constant such that $\sigma =\varepsilon $ if the
boundary condition is (\ref{bcep}) and
\[
\sigma =\varepsilon \sup_{t\in \left( 0,\delta \right) }\left| \frac{\lambda
_{i}^{L}\left( 0,t\right) }{\lambda _{i}^{R}\left( 0,t\right) }\right| <1
\]
if the boundary condition is (\ref{bceq}). Replacing $M$ by $\left( 1-\sigma
\right) M$, we can write
\[
R_{i}^{L}\left( t\right) +\hat{S}_{i}^{L}\left( t\right) \leq M\left( \left|
s_{i}^{I}\right| _{C\left[ 0,1\right] }+\left| \xi _{i}^{B}\right| _{C\left[
0,\delta \right] }+\int_{0}^{t}\left( R_{i}^{L}\left( \tau \right) +\hat{S}%
_{i}^{L}\left( \tau \right) +T_{i}^{L}\left( \tau \right) \right) d\tau
\right) .
\]
Hence, by Gronwall's inequality
\[
R_{i}^{L}\left( t\right) +\hat{S}_{i}^{L}\left( t\right) \leq M\left( \left|
s_{i}^{I}\right| _{C\left[ 0,1\right] }+\left| \xi _{i}^{B}\right| _{C\left[
0,\delta \right] }+\delta \max_{t\in \left( 0,\delta \right)
}T_{i}^{L}\left( t\right) \right) .
\]
This proves that both $R_{i}^{L}\left( t\right) $ and $S_{i}^{L}\left(
t\right) $ have upper bounds in the form of the right hand side of (\ref
{conth1.1}).

We next extend the estimate to $D_{i,\delta }^{L}$ or $D_{i,\delta }^{R}$ if
the end is either a branching junction or a transitional junction. In either
case, the solutions on the branches $j_{1},\ldots ,j_{\mu }$ connecting to
the junction constitute a fixed point of the operator $K$, which is defined
in (\ref{linth1.8}). Let
\[
W\left( t\right) =\sum_{l=1}^{\nu }\left( \hat{R}_{j_{l}}^{R}\left( t\right)
+S_{j_{l}}^{R}\left( t\right) \right) +\sum_{l^{\prime }=\nu +1}^{\mu
}\left( R_{j_{l^{\prime }}}^{L}\left( t\right) +\hat{S}_{j_{l^{\prime
}}}^{L}\left( t\right) \right)
\]
where $\hat{R}_{i}^{R}$ and $S_{i}^{R}$ are defined as in (\ref{conth1.3})
with $D_{\delta }^{C}$ substituted by $D_{\delta }^{C}\cup D_{\delta }^{R}$.
Then, from $w=Kw$ and in view of the assumption that $A_{i}$ has a positive
lower bound for all $i$ and $t>0$, we can deduce
\begin{eqnarray*}
W\left( t\right)  &\leq &\sigma \left( \sum_{l=1}^{\nu }\hat{R}%
_{j_{l}}^{R}\left( t\right) +\sum_{l^{\prime }=\nu +1}^{\mu }\hat{S}%
_{j_{l^{\prime }}}^{L}\left( t\right) \right)  \\
&&+M\left( \sum_{l=1}^{\nu }\left| r_{j_{l}}^{I}\right| _{C\left[ 0,1\right]
}+\sum_{l^{\prime }=\nu }^{\mu }\left| s_{j_{l^{\prime }}}^{I}\right|
_{C\left[ 0,1\right] }+\int_{0}^{t}\left( W\left( \tau \right) +T\left( \tau
\right) \right) d\tau \right) ,
\end{eqnarray*}
where
\[
T\left( \tau \right) =\sum_{l=1}^{\nu }T_{j_{l}}^{R}\left( \tau \right)
+\sum_{l^{\prime }=\nu +1}^{\mu }T_{j_{l^{\prime }}}^{L}\left( \tau \right)
+\sum_{l=1}^{\mu }\left| C^{j_{l}}\left( \tau \right) -\tilde{C}%
^{j_{l}}\left( \tau \right) \right|
\]
and $T_{i}^{R}\left( t\right) $ is defined as in (\ref{conth1.4}) with $%
D_{\delta }^{C}$ substituted by $D_{\delta }^{C}\cup D_{\delta }^{R}$.
Replacing $M$ by $\left( 1-\sigma \right) M$, we obtain
\[
W\left( t\right) \leq M\left( \sum_{l=1}^{\nu }\left| r_{j_{l}}^{I}\right|
_{C\left[ 0,1\right] }+\sum_{l^{\prime }=\nu }^{\mu }\left| s_{j_{l^{\prime
}}}^{I}\right| _{C\left[ 0,1\right] }+\int_{0}^{t}\left( W\left( \tau
\right) +T\left( \tau \right) \right) d\tau \right) .
\]
Hence, by Gronwall's inequality,
\[
W\left( t\right) \leq M\left( \sum_{l=1}^{\nu }\left| r_{j_{l}}^{I}\right|
_{C\left[ 0,1\right] }+\sum_{l^{\prime }=\nu }^{\mu }\left| s_{j_{l^{\prime
}}}^{I}\right| _{C\left[ 0,1\right] }+\delta \max_{t\in \left( 0,\delta
\right) }T\left( t\right) \right) .
\]
This leads to an upper bound in the form of the right hand side of (\ref
{conth1.1}) for $R_{i}^{R}\left( t\right) $, $S_{i}^{R}\left( t\right) $, $%
i=j_{1},\ldots ,j_{\nu }$, and $R_{i}^{L}\left( t\right) $, $S_{i}^{L}\left(
t\right) $, $i=j_{\nu +1},\ldots ,j_{\mu }$.

We have thus obtained an upper bound in the form of the right hand side of (%
\ref{conth1.1}) for the quantities $\left| r_{i}-\tilde{r}_{i}\right|
_{C\left( D_{\delta }\right) }$ and $\left| s_{i}-s_{i}\right| _{C\left(
D_{\delta }\right) }$. The conclusion of the lemma follows now from (\ref
{linth1.1}). \endproof

\section{The quasilinear system}

\setcounter{equation}{0} \setcounter{theorem}{0}In this section, we study
the quasilinear system where the coefficients $a_{i}$, $b_{i}$, $c_{i}$, $%
f_{i}$, $g_{i}$, $A_{i}$ and $C_{i}$ depend on both $\left( x,t\right) $ and
$\left( P_{i},Q_{i}\right) $. Under certain conditions, we show that the
system has a unique local solution. We then present a theorem on the
continuity of dependence of the solution on initial, boundary and forcing
function.

The basic idea in the proof of the existence of solution is to construct an
iterative sequence. Substituting any vector function $\left(
p_{i},q_{i}\right) $ for $\left( P_{i},Q_{i}\right) $ in $a_{i}$, etc., the
system becomes linear. Thus, we can use Theorem \ref{linth1} to get a
solution $\left( P_{i},Q_{i}\right) $. This defines a mapping $S$ from $%
u=:\left( p_{i},q_{i}\right) $ to $U=:\left( P_{i},Q_{i}\right) $, and the
solution for the quasilinear system is a fixed point of $S$. If there is a
subset of a Banach space that is invariant under $S$, then, we can construct
a sequence
\[
u_{k+1}=Su_{k},\quad k=0,1,\ldots .
\]
In the case where the limit exists and is unique, it gives rise to fixed
point of $S$. This is our approach in this section.

In this approach, conditions (\ref{cond2}) and (\ref{cond3}) are repeatedly
used. One might want to impose them for all the values of the variables.
This would give the existence and uniqueness for the global solution, as in
the case of the linear system. However, such a requirement is so restrictive
that even the original system (\ref{deaq}) cannot meet it. Therefore, we
will impose them only for $t=0$, and obtain the local solution for the
quasilinear system.

\begin{theorem}
\label{quath1}Assume that the initial and boundary functions $P_{i}^{I}$, $%
Q_{i}^{I}$, $P_{i}^{B}$, $Q_{i}^{B}$ and the system functions $a_{i}$, $b_{i}
$, $c_{i}$, $f_{i}$, $g_{i}$, $A_{i}$ and $C_{i}$ all have continuous
first-order derivatives with respect to each variable. Suppose that $a_{i}$
is positive and $A_{i}$ has a positive lower bound for all the values of
their arguments, and that conditions (\ref{cond2})--(\ref{cond3}) hold at $%
t=0$. Suppose also that the initial functions $P_{i}^{I}$, $Q_{i}^{I}$
satisfy any relevant boundary conditions at $t=0$. Then, for some $\delta >0$%
, there is a unique solution for $0\leq t<\delta $ to the quasilinear system
(\ref{depq}) with the initial and boundary conditions given by (\ref{ic})--(%
\ref{bcbq}), (\ref{bcbp'}), and (\ref{bctp})--(\ref{qc}).
\end{theorem}

\paragraph{\noindent Proof.}

We first consider the simpler case where $U^{I}=:\left( P^{I},Q^{I}\right) =0
$. Let $v=\left\{ v_{i}\right\} $, $v_{i}=\left( p_{i},q_{i}\right) $ be a
family of vector functions (not necessarily constitutes a solution) that
satisfy the initial and boundary conditions. Substitute $v$ for $U$ in the
functions $a_{i}$, $b_{i}$, $c_{i}$, $f_{i}$, $g_{i}$, $A_{i}$ and $C_{i}$.
Then, the system becomes linear and we can invoke Theorem \ref{linth1} to
obtain a solution $U$ to the linear system. This defines a mapping $S:$ $%
v\mapsto U$. A solution to the quasilinear system is then a fixed point of $S
$. We will choose a subset $X_{\delta ,M_{0}}$ of a Banach space such that
(1) $SX_{\delta ,M_{0}}\subset X_{\delta ,M_{0}}$, and (2) $S$ is
contracting in $X_{\delta ,M_{0}}$. For any scalar or vector function $f\in
C^{k}\left( D_{\delta }\right) $, let $\left| f\right| _{k,\delta }$ denote
the maximum norm of all the $k$-th order derivatives of $f$ in $D_{\delta }$%
. (If $f$ is a vector function, $\left| f\right| _{k,\delta
}=\max_{i}\left\{ \left| f_{i}\right| _{k,\delta }\right\} $.) Let $%
C_{B}\left( D_{\delta },\Bbb{R}^{2n}\right) $ denote the subset of the
vector-valued functions in $C\left( D_{\delta },\Bbb{R}^{2n}\right) $ that
satisfy the initial and boundary conditions. We seek $X_{\delta ,M_{0}}$ in
the form
\begin{equation}
X_{\delta ,M_{0}}=\left\{ v\in C_{B}\left( D_{\delta },\Bbb{R}^{2n}\right) :%
\text{ }\left| v\right| _{0,\delta }\leq M_{0},\left| v\right| _{1,\delta
}\leq M_{1}\right\}   \label{quath1.5}
\end{equation}
where $M_{0}$ is an arbitrary positive constant and $M_{1}$ is a constant to
be determined. Note that by the vanishing initial condition, for any $M_{1}$%
, $\left| U\right| _{1,\delta }\leq M_{1}$ implies $\left| U\right|
_{0,\delta }\leq M_{1}\delta $. Hence, for any $M_{0}$, we can ensure $%
\left| U\right| _{0,\delta }\leq M_{0}$ by reducing $\delta $. It remains,
therefore, only to show that for $M_{1}$ sufficiently large and $\delta $
sufficiently small, $\left| v\right| _{1,\delta }\leq M_{1}$ implies $\left|
Sv\right| _{1,\delta }\leq M_{1}$. Throughout this proof, we use $M$ to
represent any positive constant that may depend on $M_{1}$ but is otherwise
independent of $v$ and $\delta $, and use $\tilde{M}$ for any constant that
is independent of $M_{1}$, $v$ and $\delta $. The values of $M$ or $\tilde{M}
$ in different occurrences need not be equal.

Let $U=Sv$ and let $r_{i}$ and $s_{i}$ be defined by (\ref{linth1.6}). On
each branch, we show that
\begin{equation}
\max \left\{ \left| \left( r_{i}\right) _{x}\right| ,\left| \left(
s_{i}\right) _{x}\right| ,\right\} \leq M_{1}  \label{quath1.11}
\end{equation}
and
\begin{equation}
\max \left\{ \left| \left( r_{i}\right) _{t}\right| ,\left| \left(
s_{i}\right) _{t}\right| \right\} \leq M_{1}  \label{quath1.12}
\end{equation}
in $D_{\delta }^{C}$, $D_{\delta }^{L}$ and $D_{\delta }^{R}$ if $M_{1}$ is
large and $\delta $ is small. (Recall that $D_{\delta }^{C}$ etc. are the
intersections $D_{i}^{C}\cap D_{\delta }$ etc., respectively.) In fact, only
(\ref{quath1.11}) needs to be shown. To see this, first observe that the
vanishing initial condition and the compatibility of the initial and
boundary conditions gives
\[
\max_{i}\left\{ \left| P_{i}^{B}\right| _{C\left[ 0,\delta \right] },\left|
Q_{i}^{B}\right| _{C\left[ 0,\delta \right] }\right\} \leq M\delta .
\]
Hence, we obtain from Lemma \ref{conth1} with $\tilde{U}=0$ that
\begin{equation}
\left| U\right| _{0,\delta }\leq M\delta .  \label{quath1.1}
\end{equation}
From (\ref{linth1.2}) and (\ref{linth1.3}), there are constants $\tilde{M}$
and $M$ such that
\begin{equation}
\begin{array}{l}
\left| \partial _{i}^{R}r_{i}\right| \leq \left| l_{i}^{R}F_{i}\right|
+\left| \partial _{i}^{R}l_{i}^{R}\right| \left| U_{i}\right| \leq \tilde{M}%
+M\delta , \\[12pt]
\left| \partial _{i}^{L}s_{i}\right| \leq \left| l_{i}^{L}F_{i}\right|
+\left| \partial _{i}^{L}l_{i}^{L}\right| \left| U_{i}\right| \leq \tilde{M}%
+M\delta
\end{array}
\label{quath1.10}
\end{equation}
for each $i=1,\ldots ,n$. Hence, (\ref{quath1.12}) follows from (\ref
{quath1.11}), (\ref{quath1.10}) and the definition of $\partial _{i}^{L}$
and $\partial _{i}^{R}$ in (\ref{linth1.5}). We also note that (\ref
{linth1.1}) and (\ref{quath1.10}) imply
\begin{equation}
\left| \partial _{i}^{R}U_{i}\right| _{0,\delta }\leq \tilde{M}+M\delta
,\quad \left| \partial _{i}^{R}U_{i}\right| _{0,\delta }\leq \tilde{M}%
+M\delta  \label{quath1.7}
\end{equation}
for all $i$. This will be used later.

We first consider the middle region $D_{\delta }^{C}$, where the solution $%
\left( r_{i},s_{i}\right) $ satisfies the integral equations (\ref{inteq1})
and (\ref{inteq3}) with $r_{i}^{I}=s_{i}^{I}=0$. Differentiating the
equations with respect to $x$, we have
\begin{equation}
\begin{array}{l}
\displaystyle\left( r_{i}\right) _{x}=\left( l_{i}^{R}\right)
_{x}U_{i}\left( x,t\right) +\int_{0}^{t}\left[ \left( l_{i}^{R}F_{i}\right)
_{x}+\left( \partial _{i}^{R}l_{i}^{R}\right) \left( U_{i}\right)
_{x}-\left( l_{i}^{R}\right) _{x}\left( \partial _{i}^{R}U_{i}\right)
\right] \left( x_{i}^{R}\right) _{x}dt, \\
\displaystyle\left( s_{i}\right) _{x}=\left( l_{i}^{L}\right)
_{x}U_{i}\left( x,t\right) +\int_{0}^{t}\left[ \left( l_{i}^{L}F_{i}\right)
_{x}+\left( \partial _{i}^{L}l_{i}^{L}\right) \left( U_{i}\right)
_{x}-\left( l_{i}^{L}\right) _{x}\left( \partial _{i}^{L}U_{i}\right)
\right] \left( x_{i}^{L}\right) _{x}dt.
\end{array}
\label{quath1.2}
\end{equation}
Here, we used an identity from \cite[p.469]{CH62}:
\begin{equation}
\begin{array}{l}
\displaystyle\frac{d}{d\xi }\int_{a}^{b}f\left( x\left( t\right) ,t\right)
Dg\left( x\left( t\right) ,t\right) dt \\
\displaystyle\hspace{1in}=f\left( x\left( b\right) ,b\right) g_{x}\left(
x\left( b\right) ,b\right) x_{\xi }\left( b\right) -f\left( x\left( a\right)
,a\right) g_{x}\left( x\left( a\right) ,a\right) x_{\xi }\left( a\right) \\
\displaystyle\hspace{1.2in}+\int_{a}^{b}\left[ f_{x}\left( x\left( t\right)
,t\right) Dg\left( x\left( t\right) ,t\right) -Df\left( x\left( t\right)
,t\right) g_{x}\left( x\left( t\right) ,t\right) \right] x_{\xi }\left(
t\right) dt
\end{array}
\label{quath1.8}
\end{equation}
where $x\left( t\right) $ is a function such that $x\left( b\right) =\xi $
and $D=\frac{\partial }{\partial t}+x^{\prime }\left( t\right) \frac{%
\partial }{\partial x}$. (Notice that $x_{\xi }\left( b\right) =1.$) Let
\begin{equation}
R_{i}^{C}\left( t\right) =\sup_{\left\{ x:\left( x,t\right) \in D_{\delta
}^{C}\right\} }\left\{ \left| \left( r_{i}\right) _{x}\left( x,t\right)
\right| \right\} ,\quad S_{i}^{C}\left( t\right) =\sup_{\left\{ x:\left(
x,t\right) \in D_{\delta }^{C}\right\} }\left\{ \left| \left( s_{i}\right)
_{x}\left( x,t\right) \right| \right\} .  \label{quath1.9}
\end{equation}
From (\ref{quath1.1}), (\ref{quath1.7}) and (\ref{quath1.2}), we derive
\[
R_{i}^{C}\left( t\right) +S_{i}^{C}\left( t\right) \leq M\delta
+M\int_{0}^{t}\left( 1+R_{i}^{C}\left( t^{\prime }\right) +S_{i}^{C}\left(
t^{\prime }\right) \right) dt^{\prime }
\]
for $t\in \left[ 0,\delta \right] $. Hence, Gronwall's inequality gives
\[
\left| \left( r_{i}\right) _{x}\right| \leq M\delta e^{M\delta },\quad
\left| \left( s_{i}\right) _{x}\right| \leq M\delta e^{M\delta }
\]
in $D_{\delta }^{C}$. This proves (\ref{quath1.11}) in $D_{\delta }^{C}$ if $%
M_{1}$ is sufficiently large and $\delta $ is sufficiently small$.$

We next consider the left and right regions $D_{\delta }^{L}$, $D_{\delta
}^{R}$ which are next to an external end. Since the two cases are similar,
we will consider the left region only. Let $\hat{s}_{i}=s_{i}/\varepsilon $
for any $\varepsilon >0$. Then, the pair $\left( r_{i},\hat{s}_{i}\right) $
satisfies the fixed point equations of either (\ref{inteq5}) or (\ref
{inteq5'}), depending on the type of the boundary condition. Differentiating
the equations with respect to $x$ and using a slightly modified version of (%
\ref{quath1.8}) where the lower limit $a$ of the integral also depends on $%
\xi $:
\[
\begin{array}{l}
\displaystyle\frac{d}{d\xi }\int_{a}^{b}f\left( x\left( t\right) ,t\right)
Dg\left( x\left( t\right) ,t\right) dt \\
\displaystyle\hspace{1in}=f\left( x\left( b\right) ,b\right) g_{x}\left(
x\left( b\right) ,b\right) x_{\xi }\left( b\right) -f\left( x\left( a\right)
,a\right) g_{x}\left( x\left( a\right) ,a\right) x_{\xi }\left( a\right)  \\%
[8pt]
\displaystyle\hspace{1.2in}-f\left( x\left( a\right) ,a\right) Dg\left(
x\left( a\right) ,a\right) a_{\xi } \\[8pt]
\displaystyle\hspace{1.2in}+\int_{a}^{b}\left[ f_{x}\left( x\left( t\right)
,t\right) Dg\left( x\left( t\right) ,t\right) -Df\left( x\left( t\right)
,t\right) g_{x}\left( x\left( t\right) ,t\right) \right] x_{\xi }\left(
t\right) dt,
\end{array}
\]
we have
\begin{equation}
\begin{array}{l}
\displaystyle\left( r_{i}\right) _{x}=\left( \zeta
_{i}-l_{i}^{R}F_{i}-\left( \partial _{i}^{R}l_{i}^{R}\right) U_{i}\right)
\left( 0,\tau \right) \tau _{x}+\left( l_{i}^{R}\right) _{x}U_{i}\left(
x,t\right) -\left( l_{i}^{R}\right) _{x}U_{i}\left( x_{i}^{R}\right)
_{x}\left( 0,\tau \right)  \\
\displaystyle\hspace{0.5in}+\int_{\tau }^{t}\left[ \left(
l_{i}^{R}F_{i}\right) _{x}+\left( \partial _{i}^{R}l_{i}^{R}\right) \left(
U_{i}\right) _{x}-\left( l_{i}^{R}\right) _{x}\left( \partial
_{i}^{R}U_{i}\right) \right] \left( x_{i}^{R}\right) _{x}dt, \\
\displaystyle\left( \hat{s}_{i}\right) _{x}=\frac{1}{\varepsilon }\left(
l_{i}^{L}\right) _{x}U_{i}\left( t,x\right) +\frac{1}{\varepsilon }%
\int_{0}^{t}\left[ \left( l_{i}^{L}F_{i}\right) _{x}+\left( \partial
_{i}^{L}l_{i}^{L}\right) \left( U_{i}\right) _{x}-\left( l_{i}^{L}\right)
_{x}\left( \partial _{i}^{L}U_{i}\right) \right] \left( x_{i}^{L}\right)
_{x}dt,
\end{array}
\label{quath1.15}
\end{equation}
where
\[
\zeta _{i}=2\left( u_{i}P_{i}^{B}\right) _{t}+\varepsilon \left( \hat{s}%
_{i}\right) _{t}
\]
if the boundary condition is given by (\ref{bcep}), and
\[
\zeta _{i}=2\left( \frac{a_{i}u_{i}}{\lambda _{i}^{R}}Q_{i}^{B}\right)
_{t}+\varepsilon \left( \frac{\lambda _{i}^{L}}{\lambda _{i}^{R}}\right) _{t}%
\hat{s}_{i}+\varepsilon \left( \frac{\lambda _{i}^{L}}{\lambda _{i}^{R}}%
\right) \left( \hat{s}_{i}\right) _{t}
\]
if the boundary condition is given by (\ref{bceq}). This equation is valid
for any $\varepsilon $. So, we may choose $\varepsilon $ so small such that
\[
\sigma =:\varepsilon \left| \lambda _{i}^{L}\tau _{x}\left( 0,t\right)
\right| \max \left\{ 1,\left| \left( \frac{\lambda _{i}^{L}\left( 0,t\right)
}{\lambda _{i}^{R}\left( 0,t\right) }\right) \right| \right\} <1,\quad t\in
\left[ 0,\delta \right] .
\]
To proceed further, we need an estimate of $\left| \tau _{x}\left(
0,t\right) \right| $. Observe that $\tau \left( x\right) $ satisfies the
equation
\[
x_{i}^{R}\left( \tau ;x,t\right) =0
\]
where $x_{i}^{R}\left( \tau ;x,t\right) $ is the solution of the initial
value problem
\[
\frac{dx_{i}^{R}}{ds}=\lambda _{i}^{R}\left( x_{i}^{R},s\right) ,\quad
x_{i}^{R}\left( t;x,t\right) =x.
\]
By differentiation,
\begin{equation}
\lambda _{i}^{R}\left( 0,\tau \left( x\right) \right) \tau _{x}+\left. \frac{%
\partial x_{i}^{R}}{\partial x}\right| _{\left( \tau \left( x\right)
;x,t\right) }=0.  \label{quath1.16}
\end{equation}
Let $w_{i}=\partial x_{i}^{R}/\partial x$. Then, $w_{i}$ is the solution to
the linear equation
\[
\frac{dw_{i}}{ds}=\left( \lambda _{i}^{R}\right) _{x}\left( x_{i}^{R}\left(
s;x,t\right) ,s\right) w_{i},\quad w_{i}\left( t\right) =1.
\]
Solving the equation,
\[
w_{i}\left( s\right) =\exp \left( \int_{t}^{s}\left( \lambda _{i}^{R}\right)
_{x}\left( x_{i}^{R}\left( s^{\prime };x,t\right) ,s^{\prime }\right)
ds^{\prime }\right) .
\]
Returning to (\ref{quath1.16}), we find
\[
\tau _{x}=\frac{-1}{\lambda _{i}^{R}\left( 0,\tau \left( x\right) \right) }%
\exp \left( \int_{t}^{\tau \left( x\right) }\left( \lambda _{i}^{R}\right)
_{x}\left( x_{i}^{R}\left( s^{\prime };x,t\right) ,s^{\prime }\right)
ds^{\prime }\right) .
\]
Observe that $0<\tau \left( x\right) <t\leq \delta $ and the integrand is
bounded. Hence,
\begin{equation}
\left| \tau _{x}\right| \leq \tilde{M}e^{M\delta }.  \label{quath1.13}
\end{equation}
This is the estimate we need. By this estimate, for any $M_{1}$, we can
choose $\delta $ small enough such that the constants $\sigma $ and $%
\varepsilon $ are independent of $M_{1}$. Let $R_{i}^{L}\left( t\right) $
and $\hat{S}_{i}^{L}\left( t\right) $ be defined as in (\ref{quath1.9})
except that $s_{i}$ is substituted by $\hat{s}_{i}$ and $D_{\delta }^{C}$ is
substituted by $D_{\delta }^{L}\cup D_{\delta }^{C}$, We derive from (\ref
{quath1.15}) and the identity
\[
\left( \hat{s}_{i}\right) _{t}=\partial _{i}^{L}\hat{s}_{i}-\lambda
_{i}^{L}\left( \hat{s}_{i}\right) _{x}
\]
that
\[
R_{i}^{L}\left( t\right) +\hat{S}_{i}^{L}\left( t\right) \leq \sigma \hat{S}%
_{i}^{L}\left( t\right) +\tilde{M}+M\delta +M\int_{0}^{t}\left(
1+R_{i}^{L}\left( t^{\prime }\right) +\hat{S}_{i}^{L}\left( t^{\prime
}\right) \right) dt^{\prime }.
\]
Replacing $M$ and $\tilde{M}$ by $M\left( 1-\sigma \right) $ and $\tilde{M}%
\left( 1-\sigma \right) $, respectively, and applying Gronwall's inequality,
we obtain
\[
R_{i}^{L}\left( t\right) +\hat{S}_{i}^{L}\left( t\right) \leq \left( \tilde{M%
}+M\delta \right) e^{M\delta }.
\]
Since $\left| s_{i}\right| \leq \left| \hat{s}_{i}\right| $, it follows that
\[
\max \left\{ \left| \left( r_{i}\right) _{x}\right| ,\left| \left(
s_{i}\right) _{x}\right| \right\} \leq \left( \tilde{M}+M\delta \right)
e^{M\delta }
\]
in $D_{\delta }^{L}\cup D_{\delta }^{C}$. This proves (\ref{quath1.11}) in $%
D_{\delta }^{L}\cup D_{\delta }^{C}$ if $M_{1}$ is large and $\delta $ is
small.

We next consider the case where the end of the branch is a branching or
transitional junction. As before, all the branches that are connected to the
same junction are considered simultaneously. Differentiating the fixed point
equation $w=Kw$ where $w$ and $Kw$ are defined in (\ref{linth1.7}) and (\ref
{linth1.8}), respectively, we obtain (\ref{quath1.15}) in $D_{\delta
}^{L}\cup D_{\delta }^{C}$ for $i=j_{\nu +1},\ldots ,j_{\mu }$ and
\begin{equation}
\begin{array}{l}
\displaystyle\left( \hat{r}_{i}\right) _{x}=\frac{1}{\varepsilon }\left(
l_{i}^{R}\right) _{x}U_{i}\left( x,t\right) +\frac{1}{\varepsilon }%
\int_{0}^{t}\left[ \left( l_{i}^{R}F_{i}\right) _{x}+\left( \partial
_{i}^{R}l_{i}^{R}\right) \left( U_{i}\right) _{x}-\left( l_{i}^{R}\right)
_{x}\left( \partial _{i}^{R}U_{i}\right) \right] \left( x_{i}^{R}\right)
_{x}dt, \\[12pt]
\displaystyle\left( s_{i}\right) _{x}=\left( \theta
_{i}-l_{i}^{L}F_{i}-\left( \partial _{i}^{L}l_{i}^{L}\right) U_{i}\right)
\left( 1,\tau \right) \tau _{x}+\left( l_{i}^{L}\right) _{x}U_{i}\left(
x,t\right) -\left( l_{i}^{L}\right) _{x}U_{i}\left( x_{i}^{L}\right)
_{x}\left( 1,\tau \right)  \\
\displaystyle\hspace{0.5in}+\int_{\tau }^{t}\left[ \left(
l_{i}^{L}F_{i}\right) _{x}+\left( \partial _{i}^{L}l_{i}^{L}\right) \left(
U_{i}\right) _{x}-\left( l_{i}^{L}\right) _{x}\left( \partial
_{i}^{L}U_{i}\right) \right] \left( x_{i}^{L}\right) _{x}dt,
\end{array}
\label{quath1.14}
\end{equation}
in $D_{\delta }^{C}\cup D_{\delta }^{R}$ for $i=j_{1},\ldots ,j_{\nu }$,
where
\[
\begin{array}{l}
\displaystyle\zeta _{i}=\varepsilon \sum_{l=1}^{\nu }\left( n_{j_{l}}^{i}%
\hat{r}_{j_{l}}\right) _{t}\left( 1,\tau \right) +\varepsilon
\sum_{l^{\prime }=\nu +1}^{\mu }\left( n_{j_{l^{\prime }}}^{i}\hat{s}%
_{j_{l^{\prime }}}\right) _{t}\left( 0,\tau \right) +H_{i}, \\
\displaystyle\theta _{i}=\varepsilon \sum_{l=1}^{\nu }\left( m_{j_{l}}^{i}%
\hat{r}_{j_{l}}\right) _{t}\left( 1,\tau \right) +\varepsilon
\sum_{l^{\prime }=\nu +1}^{\mu }\left( m_{j_{l^{\prime }}}^{i}\hat{s}%
_{j_{l^{\prime }}}\right) _{t}\left( 0,\tau \right) +H_{i},
\end{array}
\]
and $m_{j}^{i}$, $n_{j}^{i}$ are defined in the proof of Theorem \ref{linth1}%
. Note that the estimate (\ref{quath1.13}) holds for $\tau _{x}$ in both (%
\ref{quath1.15}) and (\ref{quath1.14}), although in the latter case, $\tau $
is the $t$-coordinate of the intersection of the left-going characteristic
curve $x_{i}^{L}$ with the vertical line $x=1$. The derivation is identical.
Hence, there is a constant $\varepsilon $, independent of $M_{1}$, such that
\begin{eqnarray*}
\varepsilon \left| \tau _{x}\right| \left( \sum_{k=1}^{\nu }\left|
m_{j_{k}}^{i}\left( t\right) \right| +\sum_{k^{\prime }=\nu +1}^{\mu }\left|
m_{j_{k^{\prime }}}^{i}\left( t\right) \right| \right)  &<&1, \\
\varepsilon \left| \tau _{x}\right| \left( \sum_{k=1}^{\nu }\left|
n_{j_{k}}^{i}\left( t\right) \right| +\sum_{k^{\prime }=\nu +1}^{\mu }\left|
n_{j_{k^{\prime }}}^{i}\left( t\right) \right| \right)  &<&1
\end{eqnarray*}
in $\left[ 0,\delta \right] $. Let $\sigma $ be the maximum of the
quantities on the left hand side of the above inequalities. Define $\hat{R}%
_{i}^{R}$, $S_{i}^{R}$, $R_{i}^{L}$ and $\hat{S}_{i}^{L}$ as in (\ref
{quath1.9}) with obvious modifications. We see that the function
\[
W\left( t\right) =\sum_{l=1}^{\nu }\left( \hat{R}_{j_{l}}^{R}\left( t\right)
+S_{j_{l}}^{R}\left( t\right) \right) +\sum_{l^{\prime }=\nu +1}^{\mu
}\left( R_{j_{l^{\prime }}}^{L}\left( t\right) +\hat{S}_{j_{l^{\prime
}}}^{L}\left( t\right) \right)
\]
satisfies the inequality
\begin{eqnarray*}
\left( 1-\sigma \right) W\left( t\right)  &\leq &\sum_{l=1}^{\nu }\left(
\left( 1-\sigma \right) \hat{R}_{j_{l}}^{R}\left( t\right)
+S_{j_{l}}^{R}\left( t\right) \right) +\sum_{l^{\prime }=\nu +1}^{\mu
}\left( R_{j_{l^{\prime }}}^{L}\left( t\right) +\left( 1-\sigma \right) \hat{%
S}_{j_{l^{\prime }}}^{L}\left( t\right) \right)  \\
&\leq &\tilde{M}+M\delta +M\int_{0}^{t}\left( 1+W\left( t^{\prime }\right)
\right) dt^{\prime }.
\end{eqnarray*}
Hence, by rescaling and using Gronwall's inequality, we achieve
\[
W\left( t\right) \leq \left( \tilde{M}+M\delta \right) e^{M\delta }.
\]
This proves that
\[
\max \left\{ \left| \left( r_{i}\right) _{x}\right| ,\left| \left(
s_{i}\right) _{x}\right| \right\} \leq M_{1}
\]
in $D_{\delta }^{R}$ for $i=j_{1},\ldots ,j_{\nu }$ and in $D_{\delta }^{L}$
for $i=j_{\nu +1},\ldots ,j_{\mu }$ if $M_{1}$ is sufficiently large and $%
\delta $ is sufficiently small. We have thus proved (\ref{quath1.11}) in
this case.

\smallskip This completes the proof of (\ref{quath1.11}) in all cases. By
choosing appropriate values of $M_{1}$ and $\delta $, we thus obtain a set $%
X_{\delta ,M_{0}}$ in the form of (\ref{quath1.5}) which is invariant under
the mapping $S$.

We now show that $S$ is a contraction in $X_{\delta ,M_{0}}$. Let $U=Sv$, $%
\tilde{U}=S\tilde{v}$ for some $v,\tilde{v}\in X_{\delta }$, and let $W=U-%
\tilde{U}$. $W$ satisfies the vanishing initial and external boundary
conditions and its differential equations takes the form of (\ref{depq})
with the coefficients
\[
a_{i}=a_{i}\left( x,t,v\right) ,\ b_{i}=b_{i}\left( x,t,v\right) ,\
c_{i}=c_{i}\left( x,t,v\right) ,\ A_{i}=A_{i}\left( x,t,v\right)
\]
the forcing functions $f_{i}$ and $g_{i}$ replaced by
\begin{equation}
\hat{f}_{i}=:f_{i}\left( x,t,v\right) -f_{i}\left( x,t,\tilde{v}\right)
+\left( a_{i}\left( x,t,v\right) -a_{i}\left( x,t,\tilde{v}\right) \right)
\frac{\partial \tilde{Q}_{i}}{\partial x}  \label{quath1.17}
\end{equation}
and
\begin{equation}
\begin{array}{r}
\displaystyle\hat{g}_{i}=:g_{i}\left( x,t,v\right) -g_{i}\left( x,t,\tilde{v}%
\right) +\left( b_{i}\left( x,t,v\right) -b_{i}\left( x,t,\tilde{v}\right)
\right) \frac{\partial \tilde{P}_{i}}{\partial x} \\[12pt]
\displaystyle+2\left( c_{i}\left( x,t,v\right) -c_{i}\left( x,t,\tilde{v}%
\right) \right) \frac{\partial \tilde{Q}_{i}}{\partial x},
\end{array}
\label{quath1.18}
\end{equation}
respectively, and the functions $C_{i}$ in (\ref{bcbp'}) replaced by
\begin{eqnarray*}
\hat{C}_{j_{l}} &=&C_{j_{l}}\left( x,t,v\right) -C_{j_{l}}\left( x,t,\tilde{v}%
\right) +\left( A_{j_{l}}\left( x,t,v\right) -A_{j_{l}}\left( x,t,\tilde{v}%
\right) \right) \left( \tilde{P}_{j_{l}}-\tilde{P}_{junc}\right) , \\
\hat{C}_{j_{l^{\prime }}} &=&C_{j_{l^{\prime }}}\left( x,t,v\right) -C_{j_{l^{\prime }}}\left( x,t,%
\tilde{v}\right) +\left( A_{j_{l^{\prime }}}\left( x,t,v\right)
-A_{j_{l^{\prime }}}\left( x,t,\tilde{v}\right) \right) \left( \tilde{P}%
_{junc}-\tilde{P}_{j_{l^{\prime }}}\right)
\end{eqnarray*}
for $l=1,\ldots ,\nu $, $l^{\prime }=\nu +1,\ldots ,\mu $. By the Lipschitz
property and the boundedness $\left| \tilde{U}\right| _{1,\delta }\leq M_{1}$%
, there is a constant $M$ such that
\[
\left| \hat{f}\right| _{0,\delta }\leq M\left| v-\tilde{v}\right| _{0,\delta
},\quad \left| \hat{g}\right| _{0,\delta }\leq M\left| v-\tilde{v}\right|
_{0,\delta },\quad \left| \hat{C}\right| _{0,\delta }\leq M\left| v-\tilde{v}%
\right| _{0,\delta }
\]
Hence, by Theorem \ref{conth1},
\[
\left| Sv-S\tilde{v}\right| _{0,\delta }\leq M\delta \left| v-\tilde{v}%
\right| _{0,\delta }.
\]
Therefore, $S$ is contracting in $X_{\delta ,M_{0}}$ if $\delta $ is
sufficiently small.

The rest is standard (cf. e.g., \cite{CH62}). Starting with a $v_{0}\in
X_{\delta ,M_{0}}$, we generate an iterative sequence $v_{k+1}=Sv_{k}$.
Clearly, each $v_{k}$ lies in $X_{\delta ,M_{0}}$ and the sequence converges
uniformly. The limit then satisfies the integral equations in the proof of
Theorem \ref{linth1}, and hence, is differentiable. Therefore, it is the
solution of the quasilinear differential equations. This proves the
existence and uniqueness of the solution when $U^{I}=0$.

If $U^{I}\neq 0$, we regard $U^{I}$ as a vector function of $x$ and $t$ and
introduce $\tilde{U}=U-U^{I}$. It follows that $\tilde{U}$ is a solution of
the quasilinear equations (\ref{depq}) with the forcing functions $\tilde{f}%
_{i}$ and $\tilde{g}_{i}$ given by
\[
\tilde{f}_{i}=f_{i}-\left( Q_{i}^{I}\right) _{x}a_{i},\quad \tilde{g}%
_{i}=g_{i}-\left( P_{i}^{I}\right) _{x}b_{i}-\left( Q_{i}^{I}\right)
_{x}2c_{i}
\]
and the boundary functions are given by
\[
\tilde{P}_{i}^{B}=P_{i}^{B}-P_{i}^{I},\ \tilde{Q}%
_{i}^{B}=Q_{i}^{B}-Q_{i}^{I},
\]
and
\[
\tilde{C}_{j_{l}}=C_{j_{l}}+A_{j_{l}}P_{j_{l}}^{I},\quad \tilde{C}%
_{j_{l^{\prime }}}=C_{j_{l^{\prime }}}-A_{j_{l^{\prime }}}P_{j_{l^{\prime
}}}^{I}
\]
for $l=1,\ldots ,\nu $, $l^{\prime }=\nu +1,\ldots ,\mu $. Since $\tilde{U}$
has the vanishing initial values, it can be uniquely solved for an interval
of $t\in \left[ 0,\delta \right] $. This gives rise to a solution $U$.
\endproof

\paragraph{\emph{Remark:}}

Examples can be constructed to show that if the condition
(\ref{cond3}) fails at $t=0$, then, the local solution need not
exist or may be not unique. In particular, if (\ref{cond3}) fails
at a source end, then, the system is under-determined, and if it
fails at a terminal end, the system is over-determined.

We give next a result for the continuity of dependence of the solution and
its derivatives on the initial, boundary and forcing functions and their
derivatives. This follows from an argument similar to the proofs of Lemma
\ref{conth1} and Theorem \ref{quath1}.

\begin{corollary}
\label{quath2}Let $U=\left( P,Q\right) $ and $\tilde{U}=\left( \tilde{P},%
\tilde{Q}\right) $ be two solutions of the quasilinear problem of Theorem
\ref{quath1}. Suppose the conditions of that theorem hold for the initial
and boundary functions of both solutions. Then, there exists a constant $M>0$%
, independent of initial, boundary and forcing functions, such that
\begin{equation}
\begin{array}{r}
\left| U-\tilde{U}\right| _{k,\delta }\leq M\left( \left| P^{I}-\tilde{P}%
^{I}\right| _{C^{k}\left[ 0,1\right] }+\left| Q^{I}-\tilde{Q}^{I}\right|
_{C^{k}\left[ 0,1\right] }+\left| P^{B}-\tilde{P}^{B}\right| _{C^{k}\left[
0,\delta \right] }+\left| Q^{B}-\tilde{Q}^{B}\right| _{C^{k}\left[ 0,\delta
\right] }\right.  \\
+\left. \delta \left| f-\tilde{f}\right| _{C^{k}\left( \overline{D_{\delta }}%
\right) }+\delta \left| g-\tilde{g}\right| _{C^{k}\left( \overline{D_{\delta
}}\right) }+\delta \left| C-\tilde{C}\right| _{C^{k}\left[ 0,\delta \right]
}\right) .
\end{array}
\label{quath2.1}
\end{equation}
for $k=0,1$.
\end{corollary}

\paragraph{Proof.}

For $k=0$, the result follows from substituting one of the solutions into
the coefficients, modifying the forcing functions by (\ref{quath1.17})--(\ref
{quath1.18}), and using Lemma \ref{conth1}. For $k=1$, we differentiate the
equations and apply the lemma to the resulting equations for the derivatives
of the solution. The process is standard and is omitted. \endproof

\end{document}